\documentclass{cocv}
\usepackage{mathbbold,bbm}
\usepackage{amsmath}
\usepackage{mathrsfs,amssymb,amsmath,booktabs,array,xcolor,url,cite,color,soul,multirow,multicol}
\usepackage{slashbox}
\usepackage{stfloats}
\usepackage{amsfonts}
\usepackage{extarrows}
\usepackage{textcomp}
\usepackage{ulem}
\usepackage{cases}
\usepackage{bm}
\usepackage{graphicx}
\usepackage{arydshln}
\usepackage[all]{xy}
\usepackage[colorlinks,linkcolor=red,anchorcolor=blue,citecolor=green]{hyperref}
\newtheorem{definition}{Definition}
\newtheorem{remark}{Remark}
\newtheorem{lemma}{Lemma}
\newtheorem{proposition}{Proposition}
\newtheorem{theorem}{Theorem}

\def\dd{\text{d}}
\begin{document}

\title{Time-domain Boundedness of Noise-to-State Exponentially Stable Systems}\thanks{This work was supported by the National Nature Science Foundation of China under Grant No. 11671418, and the Zhejiang Provincial Natural Science Foundation of China under Grant No. LZ20A010002.}

\author{Zhou Fang$^1$}
\author{Chuanhou Gao}
\address{School of Mathematical Sciences, Zhejiang University,
	Hangzhou 310027, P. R. China. (\email{zhou\_fang@zju.edu.cn,
		gaochou@zju.edu.cn (Correspondence)}).}

\date{}
\begin{abstract}
In this paper we prove the time-domain boundedness for noise-to-state exponentially stable systems, and further make an estimation of its lower bound function, which allows to answer the question that how long the solution of a stochastic noise-to-state exponentially stable system stays in the domain of attraction and what happens with it if it escapes from this region for a while. The results will complement the probability-domain boundedness of noise-to-state exponentially stable systems, and provide a new insight into noise-to-state exponential stability.
\end{abstract}
\maketitle
\section{Introduction}
{Stochastic phenomena arise abundantly in real systems ranging from big power grid, chemical reaction process to financial market. For these stochastic systems, a universal phenomenon is that the noise term does not vanish or decay with time on the set of equilibria of the underlying deterministic system. The persistent noise causes a great technical obstacle for applying some classic stochastic stabilization theories, such as stochastic Lyapunov theorem\mbox{\cite{kushner67}}, stochastic passivity theorem\mbox{\cite{florchinger1999passive}}, etc., to this class of stochastic systems. It thus becomes difficult to capture stochastic stability, stochastic asymptotic stability or $p$-th moment asymptotic stability on the above set of equilibria. Even worse, the persistent noise makes it impossible at all for these stochastic systems to attain stability in the sense of probability\mbox{\cite{fang2016stochastic}}.}

To deal with persistent noise, Deng et al. proposed the notion of noise-to-state stability (NSS)\mbox{\cite{deng2001stabilization}} which concerns the ultimate boundedness instead of the stochastic convergence of state. NSS is a stochastic counterpart of the deterministic input-to-state stability\mbox{\cite{sontag1989smooth}}, where the increment of the Lyapunov function is composed of a radially unbounded dissipative term and a positive term depending on the disturbance variance. As a result, if the disturbance variance is bounded, then the state is bounded in probability, which characterizes a certain sense of stability. The innovative effect has triggered an intensive interest in the study of NSS. In the theoretical aspect, Mateos-N$\acute{\text{u}}\tilde{\text{n}}$ez and Cort$\acute{\text{e}}$s\mbox{\cite{mateos2014p}} extended the current NSS in $1$st moment\mbox{\cite{deng2001stabilization}} for stochastic systems with persistent noise to NSS in $p$th moment by proposing the concept of the $p$th NSS-Lyapunov function. In the application aspect, Ferreira et al.\mbox{\cite{ferreira2012stability}} proposed an adapted version of NSS/NSES, based on which they gave a genius criterion to say interconnected stochastic systems noise-to-state stable, and then applied it to large-scale biological reaction networks. Mateos-N$\acute{\text{u}}\tilde{\text{n}}$ez and Cort$\acute{\text{e}}$s\mbox{\cite{mateos2016noise}} found the NSS behavior in a distributed convex optimization algorithm for multi-agents systems and innovatively discussed the stability of algorithms in the perspective of NSS. Facing complicated random signals and switch behaviors in practical systems, Zhang et al.\mbox{\cite{zhang2016noise}} extended the application of NSS to random switched systems where disturbances are not necessarily white noises, and witnessed success in controlling stochastic mechanical systems. Except for NSS, there is another stability analysis strategy to accommodate the persistent noise in stochastic systems, which requires the dissipation of the Lyapunov function only outside a neighborhood of the equilibrium\mbox{\cite{zakai1969lyapunov,satoh2014bounded,fang2016stochastic}}.
{This alternative strategy can also suggest some good properties, such as stability of stationary probability distribution\mbox{\cite{zakai1969lyapunov}}, bounded stability of the state in probability\mbox{\cite{satoh2014bounded}}, the convergence of the transition measure of the state, and the ergodicity\mbox{\cite{fang2016stochastic}}.} In fact, according to \textit{Proposition 2.5} in \cite{ferreira2012stability}, NSS can be implied by the alternative method, and therefore these two analysis techniques are quite connected.

Although NSS provides an insight on the confidence level that the state stays within a bounded set, the related studies\mbox{\cite{deng2001stabilization,mateos2014p}} mainly focused on the point-wise analysis, i.e., the state boundedness in probability is estimated at a fixed time. Little information is known about how long the state will stay in the bounded set, and what happens with it if it escapes from this set. After all, from the notion of NSS there exist some risks that the state does not lie in the bounded set. These puzzles motivate us to analyze the path-wise behavior of noise-to-state stable systems, i.e., exploring {the ratio of resident time in a region.} {Particularly, if the ratio of resident time in a bounded region is lower bounded by a big number, then it indicates a particular kind of stability for NSS systems, i.e., the trajectory will evolve in this region with a large proportion of time.}
Thus, the path-wise analysis can solve the above puzzles.

It should be noted that the analysis for the {resident time ratio} is not difficult if the noise is nonsingular. Under this condition, the stochastic process will admit the Feller and irreducible properties\mbox{\cite{stettner19943}}. These properties together with the finite time recurrence of noise-to-state stable systems imply ergodicity \cite{khas1960ergodic}, which demonstrates the time average equaling to the spacial (probability) average.
In other words, the ratio of resident time in a region (or, equivalently, the average time to stay in this region) equals to the stationary probability measure of this region.
However, the analysis for general noise-to-state stable systems which may admit singular noises is still uncovered.

For the above reasons, this paper contributes to analyzing the resident time ratio for general noise-to-state exponentially stable systems, a kind of special noise-to-state stable systems that have noise-to-state exponential stability (NSES). To be specific, we intend to find a lower bound for this ratio and, therefore, provide a particular kind of stability for noise-to-state exponentially stable systems. Through borrowing loop techniques \cite{khas1960ergodic}, we replace the non-singularity condition with an exponential dissipation condition to make the current analysis scheme applicable to general noise-to-state exponentially stable systems. In this work, we firstly construct loops for the trajectory of the noise-to-state exponentially stable system and then divide the time consumed in each loop into the up crossing time and the down crossing time. Further, based on the structure of NSES, the probability distributions of the up crossing time and down crossing time are proved to be uniformly controlled by two distribution functions, respectively, which immediately suggests the estimations of the averages of the up crossing time and down crossing time.
Finally, the lower bound of the ratio of resident time is shown through the up crossing and down crossing time. Also, we show that this lower bound closely depends on the noise intensity.
The lower bound will go to $1$ as the noise intensity tends to $0$.
Remarkably, the lower boundedness of resident time ratio is an almost surely result, which we can trust without taking any risk. Therefore, from the viewpoint of probability domain, although there exist some risks for the state to escape from a bounded set, the boundedness of resident time ratio tells that the state will come back to the set in the long term if it gets away from the bounded region. It thus depicts a stable behavior of stochastic systems with non-vanishing noise in a new perspective.
In general, the analysis of resident time ratio on NSES provides a good supplement to the present probability-domain analysis\mbox{\cite{deng2001stabilization}} in revealing the essential phenomenon of NSES.

The rest of this paper is organized as follows.
Section \ref{Section NSS} briefly reviews the basic concepts and results about NSS/NSES, and then formulates the problem in study.
In Section {\ref{Section Loops}, we define loops for segmenting the trajectory along the timeline and show the probability distributions of the up/down crossing time to be uniformly controlled. Section \ref{Section Boundedness analysis} analyzes the lower boundedness of resident time ratio by estimating its lower bound all through evolution. In Section \ref{Section numerical example}, a numerical example is presented to illustrate our work. Finally, Section \ref{Section conclusion} contains the conclusions and topics of future research.

Here are some notations that readers will find in the context.

\noindent
\begin{tabular}{ll}
	{\sc Notations}: &\\
	\hline
	$\mathcal{K}$: &
	$\mathcal{K}=\left\{\rho(\cdot)\in\mathcal{C}(\mathbb{R})~|~\rho(\cdot)~\text{is strictly increasing and}~ \rho(0)=0\right\}$.\\
	$\mathcal{K}_{\infty}$: &
	$\mathcal{K}_{\infty}=\left\{\rho(\cdot)\in\mathcal{K}~|~\lim_{r\to+\infty}\rho(r)=+\infty\right\}$.\\
$\mathcal{KL}$:& $\mathcal{KL}=\{\rho(\cdot,\cdot)\in\mathcal{C}(\mathbb{R}_{\geq 0}\times \mathbb{R}_{\geq 0};\mathbb{R}_{\geq 0})~|~\text{for each fixed } s\geq 0,~\rho(\cdot,s)\in \mathcal{K},$ \\& $\text{and, for each fixed } r\geq 0, \rho(r,\cdot) \text{~is decreasing and~}  \lim_{s\to+\infty}\rho(r,s)=0\}$.\\
	$\mathbbold{1}_{\{\cdot\}}$: &
	$\mathbbold{1}_{\{\cdot\}}$ is called indicator function which is valued 1 if event $\{\cdot\}$\\
	&happens, otherwise 0. \\
	$\mathcal{W}_{-1}:$ &
	Lambert W function in the lower branch. For any real number\\ & $a\in[-\text{e}^{-1},0)$, there is  $a=\mathcal{W}_{-1}(a)\text{e}^{\mathcal{W}_{-1}(a)}$ and $\mathcal{W}_{-1}(a)\leq-1$.\\
	$U(0,1):$ &
	The uniform distribution with lower bound $0$ and upper bound $1$. \\
	$\exp\{\cdot\}$ &
	The exponential function. \\
	{$\|\cdot \|_{2}$ and $\|\cdot \|_F:$} &
	{The Euclidean norm and Forbenius norm, respectively.} \\
	\hline
\end{tabular}
\section{NSS/NSES and Problem formulation}{\label{Section NSS}}
This section reviews {some basic concepts and results about NSS/NSES \cite{deng2001stabilization,ferreira2012stability}}, based on which the problem under consideration is formulated.

First of all, let us review some basic notations and concepts of the probability space \cite{mao2007stochastic}.
In this paper, we denote the probability space as $\left(\Omega, \mathcal{F},\mathrm{P}\right)$, where $\Omega$ is the event space, $\mathcal{F}$ is a $\sigma$-algebra, and $\mathrm{P}$ is the probability measure defined on the $\sigma$-algebra.
The filtration $\{\mathcal{F}_{t}\}_{t\geq 0}$ is a family of sub-$\sigma$-algebra of $\mathcal{F}$ satisfying $\mathcal{F}_{t} \subset \mathcal{F}_s \subset \mathcal{F}$ for any $0\leq t<s<\infty$.
Moreover, it is also right continuous, i.e., {$\mathcal{F}_{t}=\cap_{s>t} \mathcal{F}_{s}$} for any $t\geq 0$, and $\mathcal{F}_0$ contains all the subsets of $\Omega$ of probability $0$.
A function $f: \Omega \times [0,+\infty) \to \mathbb{R}^{n}$, where $n$ is an arbitrary positive integer, is called an $\{\mathcal{F}_{t}\}$-adapted process, if $f(\cdot,t): \Omega\to\mathbb{R}^{n}$ is $\mathcal{F}_{t}$-measurable for any $t\geq 0$.
Such $\{\mathcal{F}_{t}\}$-adapted processes includes but is not limited to the standard $n$-dimensional Wiener process, denoted as $\mathcal{B}(\cdot, \cdot)$, and the solution (if exists) of a stochastic differential equation, denoted as $x(\cdot, \cdot)$.
For a given pair $(\omega,t)\in \Omega\times [0,+\infty)$, the $\{\mathcal{F}_{t}\}$-adapted process $f(\cdot,\cdot)$ takes the value of $f(\omega,t)$.
In this paper, we also denote this value as $f(t)$ for simplicity.
\subsection{NSS/NSES}
Consider the following $n$-dimensional stochastic differential equation
\begin{equation}{\label{Eq. stochastic systems}}
 \dd x(t) = f(x)\dd t+ h(x)\Sigma(t)\dd \mathcal{B}(t)
\end{equation}
where $x(t)\in\mathbb{R}^{n}$ is the state vector at time $t$ for $t\in [0,+\infty)$, the Borel measurable functions $f:\mathbb{R}^n\to \mathbb{R}^n$ (sometimes called the drift) and $h:\mathbb{R}^n\to\mathbb{R}^{n\times q}$ (sometimes called the diffusion) are locally bounded and locally Lipschitz continuous, the matrix-valued function $\Sigma:[0,+\infty)\to \mathbb{R}^{q\times m}$ is also Borel measurable and bounded that modulates the covariance of the noise, and $\mathcal{B}(t)\in\mathbb{R}^{m}$ represents a standard Wiener process. The conditions on $f(x)$, $h(x)$ and $\Sigma(t)$ serve for guaranteeing the local existence and uniqueness of solutions for Eq. (\ref{Eq. stochastic systems}) with respect to any initial condition $x(0)=x_0$, {where $x_{0}$ is predetermined.} If $\Sigma(t)$ is known, it is usually assumed to be the identity matrix $I$ in the model, since $h(x)$ can be redefined to incorporate unchanged $\Sigma(t)$.

For simplicity, assume that $x\equiv\mathbbold{0}_n$ is a solution of the underlying deterministic dynamics $\dot{x}(t)=f(x)$, which means $f(\mathbbold{0}_n)=\mathbbold{0}_n$. If $h(\mathbbold{0}_n)=\mathbbold{0}_{n\times q}$, then the system (\ref{Eq. stochastic systems}) is of vanishing noise, otherwise, of non-vanishing noise. For the latter, many classical stochastic notions of stability, like stochastic Lyapunov theorem and stochastic LaSalle theorem, fail to work. The concept of NSS/NSES is thus proposed, instead, to the stability of the stochastic system (\ref{Eq. stochastic systems}) with non-vanishing noise as well as unknown noise intensity $\Sigma(t)$. The specific definition of NSS/NSES is as follows \cite{deng2001stabilization,ferreira2012stability}.

\begin{definition}[NSS/NSES]{\label{Def NNS}}
For a stochastic system (\ref{Eq. stochastic systems}), suppose there exists a $\mathcal{C}^{2}$ function $V:\mathbb{R}^{n}\to\mathbb{R}_{\geq 0}$ and class $\mathcal{K}_{\infty}$ functions $\alpha_{1}$, $\alpha_{2}$, $\alpha_{3}$ and $\gamma$ such that
\begin{equation}{\label{Eq. control of V}}
\alpha_{1}\left(\|x\|_{2}\right)\leq V(x) \leq \alpha_{2}\left(\|x\|_{2}\right)
\end{equation}
and
\begin{eqnarray}{\label{Eq. def NSS}}
\mathcal{L}\left[V(x)\right]
&\triangleq& \left(\frac{\partial V}{\partial x}\right)^\top f(x)
     +\frac{1}{2} \mathrm{Tr}\left\{\Sigma(t)^{\top}h(x)^{\top}\frac{\partial^{2} V}{\partial x^{2}}h(x)\Sigma(t)\right\}\notag\\
&\leq& -\alpha_{3}(\|x\|_{2})+\gamma\left(\|\Sigma(t)\|_{F}\right),
\end{eqnarray}
then the system is said to be noise-to-state stable and $V(x)$ is called the noise-to-state Lyapunov function. Particularly, if $\alpha_{3}(\|x\|_2)\geq cV(x)$ where $c\in\mathbb{R}_{\textgreater 0}$ is a positive constant, then the system is said to be noise-to-state exponentially stable.
\end{definition}

The above definition suggests that for a noise-to-state stable system, the underlying Lyapunov function $V(x)$ is radially unbounded while $\mathcal{L}V(x)$ is upper bounded by a constant $$\gamma_{\text{max}}\triangleq\gamma\left(\sup_{t}\|\Sigma(t)\|_{F}\right).$$
These properties guarantee that the explosion time of the solution of the stochastic differential equation (\ref{Eq. stochastic systems}) is infinite with each $x_0$ \cite{Narita82}. In the subsequent investigation we focus our attention on noise-to-state exponentially stable systems. An important result for the noise-to-state exponentially stable system is that its state is bounded in probability \cite{deng2001stabilization}. Namely, for $\forall~\epsilon\textgreater 0$ and $\forall~t\geq 0$, there exist a $\mathcal{KL}$ function $\beta$ and a $\mathcal{K}_{\infty}$ function $\delta$ such that
\begin{equation}\label{eq. probability nss}
\mathrm{P}\left\{\|x(t)\|_{2}<\beta(\|x_0\|_{2},t)+\delta\left(\sup_{s}\|\Sigma(s)\|_{F}\right)\right\}
\geq 1-\epsilon.
\end{equation}
This result means that at any given time $t$ the state starting from $x_0$ will keep in a bounded region with a large probability whose size depends on the initial state and the noise covariance. It can be also read that at any $t$ and for any open ball $\mathbb{B}^n(r)$ of radius $r$ centered at the origin, the state $x(t)$ of the noise-to-state exponentially stable system lies within $\mathbb{B}^n(r)$ in probability higher than $p(t,r)$, a function valued by $t$ and $r$. Clearly, NSES characterizes a certain sense of stability for each state $x(t)$, but there still remains a certain probability risk that the state $x(t)$ will be out of the open ball. Therefore, the current annotation of NSES is not enough.
{Then, we are interested in the behavior of the whole trajectory, such as how much time it takes in staying this region and what happens if the trajectory escapes from it.}

\subsection{Problem formulation}
For the above reasons, we focus on the ratio of resident time in a closed ball, $\bar{\mathbb{B}}^n(r)\triangleq\{x\in\mathbb{R}^{n}~|~\|x\|_{2}\leq r\}$.
This ratio depicts the proportion of time for the trajectory to evolve in this ball.
Mathematically, it can be evaluated by the time average of the indicator function $\mathbbold{1}_{\{\|x(t)\|_{2} \leq r\}}$, i.e., the limit behavior of $1/T \int_{0}^{T} \mathbbold{1}_{\{\|x(t)\|_{2} \leq r\}} \dd t$ as $T$ goes to infinity.
Note that the limit of this term may not exist, so we consider the inferior limit instead and define resident ratio function by
\begin{equation}{\label{eq definition of time measured distribution}}
D(r)\triangleq\liminf_{T\to+\infty} \frac{1}{T}\int_{0}^{T} \mathbbold{1}_{\{\|x(t)\|_{2} {\leq} r\}} \dd t.
\end{equation}
{Obviously, for a fixed trajectory $D(r)$ is nondecreasing, right continuous and upper bounded by one, three basic features for a distribution function.
However, as $D(r)$ is a random function depending on the trajectory, we also call it a distribution-like function.

{Since the physical significance of $D(r)$ is to estimate the mean residence time proportion of a trajectory in the closed ball $\bar{\mathbb{B}}^n(r)$ along the whole time horizon, $1-D(r)$ measures the mean residence time proportion of this trajectory out of $\bar{\mathbb{B}}^n(r)$. If $D(r)$ is large enough compared to $1-D(r)$, then during most of the time the trajectory evolutes in the closed ball $\bar{\mathbb{B}}^n(r)$ while stays out of $\bar{\mathbb{B}}^n(r)$ for very little time. At this point, even if the trajectory sometimes escapes from the ball, it is expected not to spend too long time reentering the ball. Therefore, by means of $D(r)$, NSES can also indicate a certain stability from the time-domain viewpoint.} The attention thus changes to find a right continuous lower bound function of $D(r)$, denoted by $b(r)$, such that
\begin{equation}{\label{Eq. goal0}}
D(r)\geq b(r) \quad \text{a. s.}\quad\quad \text{and} \quad\quad \lim_{r\to+\infty} b(r)=1.
\end{equation}
{A notable advantage of $b(r)$ is that it is a deterministic variable. Moreover, it is expected to change only over the ball radius $r$, but independent of the trajectory of the noise-to-state exponentially stable system. By means of $b(r)$, it is possible to estimate the proportion of time spent by any trajectory staying in the closed ball $\bar{\mathbb{B}}^n(r)$.} If $b(r)$ exists, it is possible to find a ``stable'' region in which the trajectories of the noise-to-state exponential system evolve with a high proportion of time.
For example, let $q_k(\sup_{t}\|\Sigma(t)\|_{F})$ (or $q_{k}$ for short) be the $k$ fractile of the distribution function $b(\cdot)$ under maximum noise intensity $\sup_{t}\|\Sigma(t)\|_{F}$, then the closed ball $\bar{\mathbb{B}}^n(q_{0.9})$ is the ``stable'' region in which more than $90\%$ time from $0$ to infinity the trajectories are bounded.
Therefore, similar to boundedness in probability, we can give another stability concept, boundedness in the time horizon as: for any $\epsilon>0$ there is a $\mathcal{K_{\infty}}$ function $\delta$ such that
	\begin{equation*}
	D\left(\delta\left(\sup_{t}\|\Sigma(t)\|_{F}\right)\right)
	\geq 1-\epsilon
    \quad\quad \text{{a.~s.}}
	\end{equation*}
where $\delta(\cdot)=q_{1-\epsilon}(\cdot)$.
The existence of $b(r)$ may render an alternative insight into NSES from the time-domain viewpoint.

{In general, our goal is to find a lower bound for the distribution-like random function $D(r)$ to illustrate a particular kind of stability for noise-to-state exponential stable systems.}
\section{Loops for segmenting the trajectory along timeline}{\label{Section Loops}}
In this section, we define loops to segment the trajectory of the noise-to-state exponential stable system along the timeline, and, moreover, we analyze that the distributions of any loop up and down crossing time are uniformly controlled.
\subsection{Loop}
The calculation of the resident time ration function (or, equivalently, the time average of the indicator function) is not difficult if the noise part $h(x)\Sigma(t)$ in Eq. (\ref{Eq. stochastic systems}) satisfies the nonsingular condition, that is, $h(x)\Sigma(t)\Sigma^\top(t)h^\top(x)$ is a nonsingular matrix.
In this case, the Feller and irreducible conditions are satisfied implying that the time average of the indicator function can be evaluated directly by the corresponding stationary probability measure under the condition of $t\to +\infty$.
However, if the nonsingular condition of $h(x)\Sigma(t)\Sigma^\top(t)h^\top(x)$ is not available, the ergodicity cannot be utilized for the purpose of calculation. In this common case, we follow the idea from \cite{khas1960ergodic} and construct loop for the calculation of $D(r)$.

\begin{definition}[Loops]{\label{Def loops}}
For a noise-to-state {exponentially} stable system given by Eq. (\ref{Eq. stochastic systems}) {with noise-to-state Lyapunov function $V(\cdot)$}, we denote $\tau_{0}=0$ and
\begin{align*}
	\tau_{2i+1} =& \sup \left\{t>\tau_{2i}~|~V\left(x(t)\right)   \leq V_{1}\right\},~~~ \text{if } \tau_{2i}<+\infty, \\
	\tau_{2i+2} =& \sup \left\{t>\tau_{2i+1}~|~V\left(x(t)\right) \geq V_{0}\right\}, ~~~\text{if } \tau_{2i+1}<+\infty,
\end{align*}
where $i=0,1,2,\cdots$, and $V_{0}$, $V_{1}$ are two positive constants satisfying $c^{-1}\gamma_{\mathrm{max}}<V_{0}<V_{1}$.
The trajectory $x(t)$ from time $\tau_{2i}$ to $\tau_{2i+2}$ is called the $(i+1)$th loop of the system.
Moreover, the time differences $\tau_{2i+1}-\tau_{2i}$ and $\tau_{2i+2}-\tau_{2i+1}$ are termed as the up crossing time and down crossing time of the $(i+1)$th loop, respectively.
\end{definition}

\begin{remark}
 The number of loops is almost surely infinite when the matrix defined by Eq. (\ref{Eq. stochastic systems}), equivalently,
  $h(x)\Sigma(t)\Sigma^\top(t)h^\top(x)$ is nonsingular, however this result may not apply to the case that $h(x)\Sigma(t)\Sigma^\top(t)h^\top(x)$ is singular \cite{khas1960ergodic}.
\end{remark}

\begin{remark}{\label{remark 2}}
{If we can preclude the case where $\lim_{i\to \infty} \tau_{i} <\infty $,} then the whole evolution time for the trajectory of the noise-to-state exponentially stable system (\ref{Eq. stochastic systems}) is composed of the durations of all loops. Moreover, for every loop, i.e., $\forall i\in\mathbb{Z}_{\geq 0}$, and for $t\in [\tau_{2i},\tau_{2i+1})$, the Lyapunov function satisfies $V(x(t))\textless V_1$. As a result, the inferior time average of $\mathbbold{1}_{\{V(x(t))\leq V_{1}\}}$ is lower controlled by the average of the up crossing time. This result also applies to the inferior time average of $\mathbbold{1}_{\{\|x\|_{2}\leq \alpha_1^{-1}(V_{1})\}}$ since $\alpha_1(\|x\|_{2})\leq V(x)\textless V_1$.
{The indefiniteness of $\lim_{i\to\infty}\tau_{i}$ will be verified in the end of Section 4.1. }
\end{remark}

{According to \textit{Remark 2}, our task then turns to estimate  the average of the loop up crossing time.
To achieve this goal, we first show distributions of up crossing and down crossing time to be uniformly controlled by two distribution functions, respectively.}
\subsection{Uniformly controlled distribution of every loop up crossing time}
To estimate the average of the loop up crossing time for a noise-to-state exponentially stable system, it is necessary to estimate the distribution of $\tau_{2i+1}-\tau_{2i}$, $\forall i\in\mathbb{Z}_{\geq 0}$, firstly. Towards this task, attention is turned to the following super-martingale-like inequality.

\begin{lemma}\label{LemmaGao1}
For a noise-to-state exponentially stable system (\ref{Eq. stochastic systems}), given {any time $t\geq s\geq 0$} and any stopping time $\tau\geq s$, there is
\begin{equation}{\label{eq. main super-martingale}}
\mathrm{E}\left\{\left.\left[V\left(x(t\wedge\tau)\right)-c^{-1}\gamma_{\mathrm{max}}\right]\text{e}^{c\cdot (t\wedge\tau-{s})}
~\right|~x(s)\right\}
\leq
V(x(s))-c^{-1}\gamma_{\mathrm{max}}.
\end{equation}
\end{lemma}
\begin{proof}
{Denote the first passage time from the closed ball $\{x~|~\|x||_2\leq r\}$ by $\tilde{\tau}_{r}=\inf\{t\geq s~|~\|x(t)||_2\geq r\}$\mbox{\cite{deng2001stabilization}}, then applying the Dynkin formula\mbox{\cite{resnick2013probability}} to the $\mathcal{C}^{2}$ function $V\left(x(t\wedge\tau\wedge\tilde{\tau}_{r})\right)\text{e}^{c\cdot(t\wedge\tau\wedge\tilde{\tau}_{r})}$ and further employing the inequality \mbox{\eqref{Eq. def NSS}} we have}
\begin{eqnarray}\label{eqGao1}
&&\text{E} \left[V\left(x(t\wedge\tau\wedge\tilde{\tau}_{r})\right)\text{e}^{c\cdot(t\wedge\tau\wedge\tilde{\tau}_{r}-s)}~|~x(s)\right]
\\
&=& V\left(x(s)\right)
+ {\text{E} \left[\int_{s}^{t\wedge\tau\wedge\tilde{\tau}_{r}} \text{e}^{c(\tilde{s}-s)} \Big(\mathcal{L}\left[V\left(x(\tilde{s})\right)\right]+cV(x(\tilde{s}))\Big)d \tilde{s}~\Big|~x(s)\right]} \notag\\
&\leq& V\left(x(s)\right)
+ {\text{E} \left[\int_{s}^{t\wedge\tau\wedge\tilde{\tau}_{r}} \text{e}^{c(\tilde{s}-s)} \gamma\left(\|\Sigma(\tilde{s})\|_{F}\right) d \tilde{s}~\Big|~x(s)\right]}\label{eq lemma 1 proof}\notag\\
&\leq& V\left(x(s)\right)
+ \text{E} \left[(\text{e}^{c\cdot(t\wedge\tau\wedge\tilde{\tau}_{r}-s)}-1)~|~x(s)\right]c^{-1}\gamma_{\mathrm{max}} \notag\\
&\leq& V\left(x(s)\right)
+ \text{E} \left[(\text{e}^{c\cdot(t\wedge\tau-s)}-1)~|~x(s)\right]c^{-1}\gamma_{\mathrm{max}}. \notag
\end{eqnarray}
{Note that the solution of the stochastic differential equation \mbox{(\ref{Eq. stochastic systems})} will be of no-explosion during finite time, i.e., $\lim_{r\to+\infty}\text{P}\{\tilde{\tau}_{r}=+\infty\}=1$. This means that  $\lim_{r\to+\infty}t\wedge\tau\wedge\tilde{\tau}_{r}=t\wedge\tau$ holds almost surely, and thus we obtain}
\begin{align*}
	\text{E} \left[V\left(x(t\wedge\tau)\right)\text{e}^{c\cdot(t\wedge\tau-s)}~|~x(s)\right]
	=\text{E} \left[\lim_{r\to+\infty}V\left(x(t\wedge\tau\wedge\tilde{\tau}_{r})\right)\text{e}^{c\cdot(t\wedge\tau\wedge\tilde{\tau}_{r}-s)}~|~x(s)\right].\notag
\end{align*}
{For the right hand term, employing {the famous Fatou Lemma}\mbox{\cite{resnick2013probability}} and also combining Eq. \mbox{(\ref{eqGao1})} yield}
\begin{eqnarray}
&&\text{E} \left[\lim_{r\to+\infty}V\left(x(t\wedge\tau\wedge\tilde{\tau}_{r})\right)\text{e}^{c\cdot(t\wedge\tau\wedge\tilde{\tau}_{r}-s)}~|~x(s)\right]
\notag\\
&\leq& \liminf_{r\to +\infty}\text{E} \left[V\left(x(t\wedge\tau\wedge\tilde{\tau}_{r})\right)\text{e}^{c\cdot(t\wedge\tau\wedge\tilde{\tau}_{r}-s)}~|~x(s)\right]
\notag\\
&\leq& V\left(x(s)\right)
+ \text{E} \left[(\text{e}^{c\cdot(t\wedge\tau-s)}-1)~|~x(s)\right]c^{-1}\gamma_{\mathrm{max}}, \notag
\end{eqnarray}
so we have
$$\text{E} \left[V\left(x(t\wedge\tau)\right)\text{e}^{c\cdot(t\wedge\tau-s)}~|~x(s)\right]\leq V\left(x(s)\right)
+\text{E} \left[(\text{e}^{c\cdot(t\wedge\tau-s)}-1)~|~x(s)\right]c^{-1}\gamma_{\mathrm{max}}.$$
By subtracting $\text{E} \left[\text{e}^{c\cdot(t\wedge\tau-s)}~|~x(s)\right]c^{-1}\gamma_{\mathrm{max}}$ from both sides, we get the inequality \mbox{(\ref{eq. main super-martingale})}.
\end{proof}

\begin{remark}
The inequality (\ref{eq. main super-martingale}) can be further simplified if the stopping time $\tau$ is set to be positive infinite. Under this condition, $t\wedge\tau=t$, then we get
	\begin{equation}\label{Ineq conseqence super-martingale}
	\mathrm{E}\Big\{\left.\left[V\left(x(t)\right)-c^{-1}\gamma_{\mathrm{max}}\right]\mathrm{e}^{ct}
	~\right|~x(s)\Big\}
	\leq
	\left[V\left(x(s)\right)-c^{-1}\gamma_{\mathrm{max}}\right]e^{cs}.
	\end{equation}
The above inequality suggests that the stochastic process $\left[V\left(x(t)\right)-c^{-1}\gamma_{\mathrm{max}}\right]\text{e}^{ct}$ is super-martingale,
so we call the inequality (\ref{eq. main super-martingale}) to be super-martingale-like.
\end{remark}

Based on the super-martingale-like inequality, it is possible to estimate the distribution of any loop up crossing time for a noise-to-state exponentially stable system.

\begin{proposition}{\label{Lem distribution of up crossing time}}
{If the stochastic system (\mbox{\ref{Eq. stochastic systems}}) admits NSES and $\mathrm{P}\{\tau_{2i}<+\infty\}>0$} holds true for some $i\in\mathbb{Z}_{> 0}$, then $\forall s \geq 0$ and $\forall x(\tau_{2i})\in\mathbb{R}^{n}$ there is
	\begin{equation}{\label{Ineq distribution of up crossing time}}
	\mathrm{P}\left\{\tau_{2i+1}-\tau_{2i}> s~|~{x(\tau_{2i})} \right\}\geq
	\frac{V_{1}-V_{0}}{V_{1}-c^{-1}\gamma_{\mathrm{max}}+c^{-1}\gamma_{\mathrm{max}}\mathrm{e}^{cs}},
	\end{equation}
where {$\{\tau_{2i+1}-\tau_{2i}\}$} represents the loop up crossing time series.
\end{proposition}

\begin{proof}
For any $t>\tau_{2i}$ with $\tau_{2i}\textless +\infty$, assume $0<\mathrm{P}\{\tau_{2i+1}\leq t~|~x(\tau_{2i})\}<1$, then by {\mbox{\eqref{eq. main super-martingale}}} we have
\begin{eqnarray}\label{PropersitionGao1}
  &&\text{E}\left\{\left.\left[V\left(x(t\wedge\tau_{2i+1})\right)-c^{-1}\gamma_{\mathrm{max}}\right]\text{e}^{c\cdot(t\wedge\tau_{2i+1}-\tau_{2i})}
  ~\right|~x(\tau_{2i})\right\} \\
  &=& \text{E}\left\{\left.\left[V\left(x(t\wedge\tau_{2i+1})\right)-c^{-1}\gamma_{\mathrm{max}}\right]\text{e}^{c\cdot(t\wedge\tau_{2i+1}-\tau_{2i})}
               ~\right|~\tau_{2i+1}\leq t,x(\tau_{2i})\right\} \notag \\
   &&\cdot\mathrm{P}\{\tau_{2i+1}\leq t~|~x(\tau_{2i})\} \notag\\
  && + \text{E}\left\{\left.\left[V\left(x(t\wedge\tau_{2i+1})\right)-c^{-1}\gamma_{\mathrm{max}}\right]\text{e}^{c\cdot(t\wedge\tau_{2i+1}-\tau_{2i})}
               ~\right|~\tau_{2i+1}> t,x(\tau_{2i})\right\}\notag\\
  &&\cdot\mathrm{P}\{\tau_{2i+1} > t~|~x(\tau_{2i})\} \notag\\
  &\geq& \text{E}\left[\left.V\left(x(\tau_{2i+1})\right)-c^{-1}\gamma_{\mathrm{max}}~\right|~\tau_{2i+1}\leq t,x(\tau_{2i})\right] \mathrm{P}\{\tau_{2i+1}\leq t~|~x(\tau_{2i})\} \notag\\
  && + \text{E}\left[\left.-c^{-1}\gamma_{\mathrm{max}}\cdot\text{e}^{c\cdot(t-\tau_{2i})}~\right|~\tau_{2i+1}> t,x(\tau_{2i})\right] \mathrm{P}\{\tau_{2i+1} > t~|~x(\tau_{2i})\} \notag\\
  &=&  \left(V_{1}-c^{-1}\gamma_{\mathrm{max}}\right)
                                \mathrm{P}\{\tau_{2i+1}\leq t~|~x(\tau_{2i})\} \notag\\
   &&-c^{-1}\gamma_{\mathrm{max}}\cdot\text{e}^{c\cdot(t-\tau_{2i})}
                                \left(1-\mathrm{P}\{\tau_{2i+1} \leq t~|~x(\tau_{2i})\}\right)
          \notag\\
  &=& \left(V_{1}-c^{-1}\gamma_{\mathrm{max}}+c^{-1}\gamma_{\mathrm{max}}\cdot\text{e}^{c\cdot(t-\tau_{2i})}\right) \mathrm{P}\{\tau_{2i+1}\leq t~|~\tau_{2i}<+\infty,~x(\tau_{2i})\}\notag\\
  &&-c^{-1}\gamma_{\mathrm{max}}\cdot\text{e}^{c\cdot(t-\tau_{2i})}. \notag
\end{eqnarray}
Note that the second equality follows immediately from the following reasons: (i) By \textit{Definition \mbox{\ref{Def loops}}}, $V\left(x(\tau_{2i+1})\right)-c^{-1}\gamma_{\mathrm{max}}$ $(=V_{1}-c^{-1}\gamma_{\mathrm{max}})$ is a deterministic number that results in the conditional expectation vanish, i.e.,
\begin{eqnarray}
&&\text{E}\left[\left.V\left(x(\tau_{2i+1})\right)-c^{-1}\gamma_{\mathrm{max}}~\right|~\tau_{2i+1}\leq t,x(\tau_{2i})\right] \mathrm{P}\{\tau_{2i+1}\leq t~|~\tau_{2i}<+\infty,~x(\tau_{2i})\}\notag\\
&=&\left(V_{1}-c^{-1}\gamma_{\mathrm{max}}\right)\mathrm{P}\{\tau_{2i+1}\leq t~|~\tau_{2i}<+\infty,~x(\tau_{2i})\}; \notag
\end{eqnarray}
{(ii) all expressions, including the equality (\mbox{\ref{Ineq distribution of up crossing time}}), are evaluated with $\tau_{2i}$ as the initial point, so $\tau_{2i}$ can be thought as a deterministic variable in this circumstance (although it is essentially a random variable as given in \textit{Definition 2}). This together with the condition of any given $t>\tau_{2i}$ indicates that $c^{-1}\gamma_{\mathrm{max}}\cdot\text{e}^{c\cdot(t-\tau_{2i})}$ is also deterministic. As a result,}
\begin{eqnarray}
&& \text{E}\left[\left.-c^{-1}\gamma_{\mathrm{max}}\cdot\text{e}^{c\cdot(t-\tau_{2i})}~\right|~\tau_{2i+1}> t,x(\tau_{2i})\right] \mathrm{P}\{\tau_{2i+1} > t~|~x(\tau_{2i})\}\notag\\
&=& -c^{-1}\gamma_{\mathrm{max}}\cdot\text{e}^{c\cdot(t-\tau_{2i})}
                                \left(1-\mathrm{P}\left\{\tau_{2i+1} \leq t~|~x(\tau_{2i})\right\}\right). \notag
\end{eqnarray}
Through applying the inequality (\ref{eq. main super-martingale}) to Eq. (\mbox{\ref{PropersitionGao1}}) with $\tau_{2i}$ as the initial time, we get
\begin{equation*}
\mathrm{P}\{\tau_{2i+1}\leq t~|~x(\tau_{2i})\} \leq
\frac{V_{0}-c^{-1}\gamma_{\mathrm{max}}+c^{-1}\gamma_{\mathrm{max}}\cdot\text{e}^{c\cdot(t-\tau_{2i})}}
{V_{1}-c^{-1}\gamma_{\mathrm{max}}+c^{-1}\gamma_{\mathrm{max}}\cdot\text{e}^{c\cdot(t-\tau_{2i})}}.
\end{equation*}
In addition, $\tau_{2i}$ is thought as a deterministic variable in the current circumstance, so we can set $s=t-\tau_{2i}$ and replace $t$ with $s+\tau_{2i}$ in the above equation, which yields the equality (\ref{Ineq distribution of up crossing time}).

Finally, we consider the special case of $\mathrm{P}\{\tau_{2i+1}\leq t~|~x(\tau_{2i})\}=0$ for any $t>\tau_{2i}$ with $\tau_{2i}\textless +\infty$. By letting $s=t-\tau_{2i}$ we get $\mathrm{P}\{\tau_{2i+1}-\tau_{2i}>s~|~x(\tau_{2i})\}=1,$ which supports the equality \mbox{(\ref{Ineq distribution of up crossing time})} without doubt. In addition, it is impossible to have $\mathrm{P}\{\tau_{2i+1}\leq t~|~x(\tau_{2i})\}=1$ for any $t>\tau_{2i}$ with $\tau_{2i}\textless +\infty$ while $\tau_{2i+1}>\tau_{2i}$ from \textit{Definition \mbox{\ref{Def loops}}}. Therefore, the result is true.
\end{proof}

Inequality (\ref{Ineq distribution of up crossing time}) gives the lower bound estimation of the distribution of the up crossing time for every loop under the condition of $\tau_{2i}\textless +\infty$ ({$i\in\mathbb{Z}_{>0}$}). In order to make sure that these distributions are uniformly controlled, we construct a random variable $\tilde{X}$ with the survival function
\begin{equation*}{\label{eq define tilde X}}
\mathrm{P}\{\tilde{X}>s\}=
\left\{
\begin{array}{cl}
\frac{V_{1}-V_{0}}{V_{1}-c^{-1}\gamma_{\mathrm{max}}+c^{-1}\gamma_{\mathrm{max}}\text{e}^{cs}}, & s\geq 0; \\
1,   &  s<0.
\end{array}	
\right.
\end{equation*}
Obviously, this survival function will decay exponentially as $s$ increases. The expectation of $\tilde{X}$ is thus finite and can be calculated by
\begin{eqnarray}{\label{eq calculate t_{uc}}}
\text{E}[\tilde{X}]
       &=&\int_{0}^{+\infty} \mathrm{P}\{\tilde{X}>s\} d s\\
        &=& \int_{0}^{+\infty}
                \frac{V_{1}-V_{0}}{V_{1}-c^{-1}\gamma_{\mathrm{max}}} \cdot
                \frac{V_{1}-c^{-1}\gamma_{\mathrm{max}}}
                {V_{1}-c^{-1}\gamma_{\mathrm{max}}+c^{-1}\gamma_{\mathrm{max}}\cdot\text{e}^{cs}} d s \notag\\
              &=& c^{-1}\frac{V_{1}-V_{0}}{V_{1}-c^{-1}\gamma_{\mathrm{max}}}\ln{\frac{V_{1}}{c^{-1}\gamma_{\mathrm{max}}}}, \notag
\end{eqnarray}
which is denoted by $t_{\text{u}}$ in the context for convenience of quotation. The finiteness of $t_{\text{u}}$ together with \textit{Proposition 1} suggests that under the condition of $\tau_{2i}\textless +\infty$, the distribution of any loop up crossing time is uniformly controlled by that of the integrable random variable, i.e., $\forall s\geq 0$ {and $\forall x(\tau_{2i})\in\mathbb{R}^{n}$},
there is
\begin{equation}\label{aa}
\mathrm{P}\left\{\tau_{2i+1}-\tau_{2i}> s~|~x(\tau_{2i}) \right\}\geq\mathrm{P}\{\tilde{X}>s\}.
\end{equation}

Although the above result is constructive, it cannot cover the case that there are finitely many loops in the noise-to-state exponentially stable system, i.e., $\tau_{2i}=+\infty$. To address this issue, we define a sequence of random variables $\{\tilde{X}_{i}\} (i\in\mathbb{Z}_{> 0}$) by
\begin{equation*}
   \tilde{X}_{i}=
   \left\{
   \begin{array}{cl}
    \tau_{2i+1}-\tau_{2i}, & \tau_{2i}<+\infty; \\
    +\infty,               & \tau_{2i}=+\infty.
   \end{array}
   \right.
\end{equation*}
which naturally represent the loop up crossing time series when the number of loops is infinite. Obviously, the conditional probability $\mathrm{P}\left\{\tilde{X}_{i}>s~|~\tau_{2i}=+\infty,\tilde{X}_{1},\cdots,\tilde{X}_{i-1}\right\}$ equals to $1$ if $s<+\infty$, otherwise $0$. Hence, for any $s\in\mathbb{R}_{\geq 0}$ and any $i\in\mathbb{Z}_{\textgreater 0}$ there is
\begin{equation}{\label{Ineq estimation of tilde Xi}}
  \mathrm{P}\left\{\tilde{X}_{i}>s~|~\tau_{2i}=+\infty,\tilde{X}_{1},\cdots,\tilde{X}_{i-1}\right\}
  \geq \mathrm{P}\left\{\tilde{X}>s\right\}.
\end{equation}
Note that
\begin{eqnarray}
  &&\mathrm{P}\left\{\tilde{X}_{i}> s~
       |~\tilde{X}_{1},\cdots,\tilde{X}_{i-1}\right\}\notag\\
         &=&  \mathrm{P}\left\{\tilde{X}_{i}> s~|~\tau_{2i}=+\infty,\tilde{X}_{1},\cdots,\tilde{X}_{i-1}\right\}
       \mathrm{P}\left\{\tau_{2i}=+\infty~
             |~\tilde{X}_{1},\cdots,\tilde{X}_{i-1}\right\}\notag\\
  && + \mathrm{P}\left\{
      \tilde{X}_{i}>  s \notag
      ~|~\tau_{2i}<+\infty,\tilde{X}_{1},\cdots,\tilde{X}_{i-1}\right\}
       \mathrm{P}\left\{\tau_{2i}<+\infty~
             |~\tilde{X}_{1},\cdots,\tilde{X}_{i-1}\right\} \notag \\
   &=&  \mathrm{P}\left\{\tilde{X}_{i}> s~|~\tau_{2i}=+\infty,\tilde{X}_{1},\cdots,\tilde{X}_{i-1}\right\}
         \mathrm{P}\left\{\tau_{2i}=+\infty~
         |~\tilde{X}_{1},\cdots,\tilde{X}_{i-1}\right\}\notag\\
     && +  \text{E}\left[\left.\mathrm{P}\left\{
       \tilde{X}_{i}>  s \notag
       ~| ~x(\tau_{2i})\right\}~\right|~ \tau_{2i}<+\infty,\tilde{X}_{1},\cdots,\tilde{X}_{i-1}\right]
      \mathrm{P}\left\{\tau_{2i}<+\infty~
       |~\tilde{X}_{1},\cdots,\tilde{X}_{i-1}\right\}, \notag
\end{eqnarray}
where, analogous to Eq. \mbox{(\ref{PropersitionGao1})}, $\mathrm{P}\left\{\tau_{2i}=+\infty~
             |~\tilde{X}_{1},\cdots,\tilde{X}_{i-1}\right\}$ is assumed to be in (0,1). Combing above equalities together with the inequalities (\ref{aa}) and (\ref{Ineq estimation of tilde Xi}) yields
\begin{equation}{\label{eq. important 1}}
	\mathrm{P}\left\{\tilde{X}_{i}> s~
	|~\tilde{X}_{1},\cdots,\tilde{X}_{i-1}\right\}\geq \mathrm{P}\left\{\tilde{X}> s\right\},
\end{equation}
i.e., the series $\{\tilde{X}_{i}\}$ ($i\in\mathbb{Z}_{> 0}$) is also controlled by $\tilde{X}$ from below.

The extreme cases that $\mathrm{P}\left\{\tau_{2i}=+\infty~|~\tilde{X}_{1},\cdots,\tilde{X}_{i-1}\right\}$ is $0$ or $1$ also support the above inequality.

{Note that there is a jump for the distribution of $\tilde{X}$ at the point 0, which at first glance might invalidate above ``control" inequality.
However, the jump behavior actually reduce the probability $\mathrm{P}\{\tilde{X}> s\}$ at $s=0$, which make things nicer in providing a lower bound.
In this sense, the jump behavior takes advantages to above analysis, and the inequality \mbox{\eqref{eq. important 1}} is not disturbed by this factor.
}
\subsection{Uniformly controlled distribution of every loop down crossing time}
The inferior time average of the indicator function also depends on the average of the loop down crossing time. Similar to the method dealing with the loop up crossing time series, we manage to analyze that the distribution of the down crossing time for every loop is uniformly controlled too.

Consider the loop down crossing time series $\{\tau_{2i}-\tau_{2i-1}\}~(i\in\mathbb{Z}_{\textgreater 0})$ for a noise-to-state exponentially stable system. It can be easily obtained that the down crossing time of any loop is almost surely finite if the number of loops is finite.

\begin{lemma}{\label{Lem down crossing time is finite}}
For a noise-to-state exponentially stable system (\ref{Eq. stochastic systems}), if for some $i\in\mathbb{Z}_{\textgreater 0}$ there is {$\mathrm{P}\{\tau_{2i-1}<+\infty\}>0$}, then the $i$-th loop down crossing time is almost surely finite, i.e.,
	\begin{equation*}
	     \mathrm{P}\{\tau_{2i}-\tau_{2i-1}<+\infty~|~\tau_{2i-1}<+\infty\}=1.
	\end{equation*}
\end{lemma}

\begin{proof}
Since the trajectory of the noise-to-state exponentially stable system (\mbox{\ref{Eq. stochastic systems}}) is continuous, for any $s\in[\tau_{2i-1},t\wedge\tau_{2i}]$ ($t>\tau_{2i-1}$) there are $V\left(x(s)\right)\geq V_{0}$ and $V(x(\tau_{2i-1}))=V_1$. {Analogous to the proof of \textit{Lemma} \mbox{\ref{LemmaGao1}}, we define a random variable $\hat{\tau}_{r}=\inf\{t\textgreater \tau_{2i-1}~|~\|x(t)||_2\geq r\}$ to represent the first passage time from the bounded set $\{x~|~\|x\|_2\leq r\}$ after $\tau_{2i-1}$. Hence, we obtain}
	\begin{eqnarray}{\label{eq lemma 2 proof}}
	  &&\text{E}\left[
	      V\left(x(t\wedge\tau_{2i}\wedge\hat{\tau}_{r})\right)~|~\tau_{2i-1}<+\infty
	      \right] \\
	   &=& \text{E}\left[V(x(\tau_{2i-1}))~|~\tau_{2i-1}<+\infty\right]
	   +{\text{E}\left[\left.\int_{\tau_{2i-1}}^{t\wedge\tau_{2i}\wedge\hat{\tau}_{r}}
	   \mathcal{L} V(x(s)) ds ~\right|~ \tau_{2i-1}<+\infty\right]} \notag\\
	   &\leq& V_{1}+\left(-cV_{0}+\gamma_{\mathrm{max}}\right)	   	   \text{E}\left[t\wedge\tau_{2i}\wedge\hat{\tau}_{r}-\tau_{2i-1}~|~\tau_{2i-1}<+\infty\right]\notag\\
       &\leq& V_{1}+\left(-cV_{0}+\gamma_{\mathrm{max}}\right)
	   	   \text{E}\left[t\wedge\tau_{2i}-\tau_{2i-1}~|~\tau_{2i-1}<+\infty\right]. \notag
	\end{eqnarray}
{Following the same proof as given in \textit{Lemma} \mbox{\ref{LemmaGao1}}, we have}
\begin{eqnarray}
 &&\text{E}\left[
	      V\left(x(t\wedge\tau_{2i}\right)~|~\tau_{2i-1}<+\infty
	      \right]\notag\\
 &=&  \text{E}\left[ \lim_{r\to\infty}
	      V\left(x(t\wedge\tau_{2i}\wedge\hat{\tau}_{r})\right)~|~\tau_{2i-1}<+\infty
	      \right]\notag\\
 &\leq& \liminf_{r\to\infty} \text{E}\left[
	      V\left(x(t\wedge\tau_{2i}\wedge\hat{\tau}_{r})\right)~|~\tau_{2i-1}<+\infty
	      \right]\notag \\
 &\leq& V_{1}+\left(-cV_{0}+\gamma_{\mathrm{max}}\right)
	   	   \text{E}\left[t\wedge\tau_{2i}-\tau_{2i-1}~|~\tau_{2i-1}<+\infty\right].\notag
\end{eqnarray}
{Note that $V_{0}>c^{-1}\gamma_{\mathrm{max}}$, which can be seen from \textit{Definition \mbox{\ref{Def loops}}}, thus we get}
	\begin{eqnarray*}
		\text{E}\left[t\wedge\tau_{2i}-\tau_{2i-1}~|~\tau_{2i-1}<+\infty\right]&\leq&
		\frac{V_{1}-\text{E}\left[
			V\left(x(t\wedge\tau_{2i})\right)~|~\tau_{2i-1}<+\infty
			\right]}
		{cV_{0}-\gamma_{\mathrm{max}}}\\
		&\leq& \frac{V_{1}}{{cV_{0}-\gamma_{\mathrm{max}}}}.
	\end{eqnarray*}
{Also, since $\lim_{t\to+\infty}t\wedge\tau_{2i}-\tau_{2i-1}=\tau_{2i}-\tau_{2i-1}$, we apply the monotone convergence theorem\mbox{\cite{resnick2013probability}} that states ``for a monotone increasing nonnegative random series, if they converges to a random variable, then the expectation of the series also converges to the expectation of this random variable" to yield}
\begin{eqnarray}
\text{E}\left[\tau_{2i}-\tau_{2i-1}~|~\tau_{2i-1}<+\infty\right]
&=&\lim_{t\to\infty}\text{E}\left[t\wedge\tau_{2i}-\tau_{2i-1}~|~\tau_{2i-1}<+\infty\right] \notag \\
&\leq& \frac{V_{1}}{{cV_{0}-\gamma_{\mathrm{max}}}} \notag \\
&<&+\infty \notag.
\end{eqnarray}
{Therefore, the $i$th ($\forall i\in\mathbb{Z}_{\textgreater 0}$) loop down crossing time, $\tau_{2i}-\tau_{2i-1}$, is almost surely finite under the condition $\tau_{2i-1}<+\infty$.}
\end{proof}

The distribution of the down crossing time of every loop is estimated as follows.
\begin{proposition}{\label{Lem distribution of down crossing time}}
For a noise-to-state exponentially stable system (\ref{Eq. stochastic systems}), if for some $i\in\mathbb{Z}_{\textgreater 0}$ there is $\mathrm{P}\{\tau_{2i-1}<+\infty\}>0$, then
 $\forall s\geq \frac{1}{c}\ln\frac{V_{1}-c^{-1}\gamma_{\mathrm{max}}}{V_{0}-c^{-1}\gamma_{\mathrm{max}}}$ {and $\forall x(\tau_{2i-1})\in\mathbb{R}^{n}$}, the $i$-th down crossing time satisfies
	\begin{equation}{\label{Ineq distribution of down crossing time}}
	   \mathrm{P}\{\tau_{2i}-\tau_{2i-1}\geq s~|~x(\tau_{2i-1})\} \leq  \frac{V_{1}-c^{-1}\gamma_{\mathrm{max}}}{V_{0}-c^{-1}\gamma_{\mathrm{max}}} \text{e}^{-cs}.
	\end{equation}
\end{proposition}

\begin{proof}
According to \textit{Lemma \ref{Lem down crossing time is finite}}, the down crossing time of any loop is almost surely finite, so we can obtain $\text{E}\left[V(x(\tau_{2i}))~|~\tau_{2i-1}<+\infty\right]=V_{0}$, $\forall i\in\mathbb{Z}_{\textgreater 0}$ and $\mathrm{P}\{\tau_{2i-1}<+\infty\}>0$. Again, we apply the inequality (\ref{eq. main super-martingale}) by setting $\tau$, the initial time and $t$ as $\tau_{2i}$, $\tau_{2i-1}$ and $+\infty$, respectively, and can get
\begin{equation*}
\text{E}\left[\left(V_{0}-c^{-1}\gamma_{\mathrm{max}}\right)\text{e}^{c\cdot(\tau_{2i}-\tau_{2i-1})}~|~x(\tau_{2i-1})\right]
		\leq V_{1}-c^{-1}\gamma_{\mathrm{max}}.
\end{equation*}
{Furthermore}, by dividing the above inequality by $V_{0}-c^{-1}\gamma_{\mathrm{max}}$,  a positive number, on both sides of the inequality, we obtain
	\begin{equation*}
		\text{E}\left[\text{e}^{c\cdot(\tau_{2i}-\tau_{2i-1})}~|~x(\tau_{2i-1})\right]
		\leq \frac{V_{1}-c^{-1}\gamma_{\mathrm{max}}}
		          {V_{0}-c^{-1}\gamma_{\mathrm{max}}}.
	\end{equation*}
Finally, utilizing {the Markov inequality} we have
	\begin{equation*}
		\mathrm{P}\left\{\text{e}^{c\cdot(\tau_{2i}-\tau_{2i-1})}\geq \text{e}^{cs}~|~x(\tau_{2i-1})\right\}
		\leq \frac{V_{1}-c^{-1}\gamma_{\mathrm{max}}}
		{V_{0}-c^{-1}\gamma_{\mathrm{max}}}\text{e}^{-cs}.
	\end{equation*}
	which implies the result of Eq. (\ref{Ineq distribution of down crossing time}) {to be} true.
\end{proof}

In the following, we adopt the same strategy as in \textbf{Sec. 3.2} to prove the distribution of the down crossing time for every loop to be uniformly controlled under the condition of $\tau_{2i-1}<+\infty$. We construct a random variable $\hat{X}$ with the distribution
\begin{equation*}{\label{eq define hat X}}
	\mathrm{P}\left\{\hat{X}>s\right\}=
	\left\{
	\begin{array}{cl}
	\frac{V_{1}-c^{-1}\gamma_{\mathrm{max}}}{V_{0}-c^{-1}\gamma_{\mathrm{max}}}\text{e}^{-cs},  &   s\geq \frac{1}{c}\ln\frac{V_{1}-c^{-1}\gamma_{\mathrm{max}}}{V_{0}-c^{-1}\gamma_{\mathrm{max}}}; \\
	1,  &  s< \frac{1}{c}\ln\frac{V_{1}-c^{-1}\gamma_{\mathrm{max}}}{V_{0}-c^{-1}\gamma_{\mathrm{max}}}.
	\end{array}
	\right.
\end{equation*}
Its expectation can be accordingly calculated as
\begin{eqnarray}{\label{eq calculate t_{dc}}}
\text{E}[\hat{X}]
      &=&\int_{0}^{+\infty} \mathrm{P}(\hat{X}>s) \dd s\\
      &=& \int_{0}^{c^{-1}\ln{\frac{V_{1}-c^{-1}\gamma_{\mathrm{max}}}{V_{0}-c^{-1}\gamma_{\mathrm{max}}}}} 1 \dd s
        + \int_{c^{-1}\ln{\frac{V_{1}-c^{-1}\gamma_{\mathrm{max}}}{V_{0}-c^{-1}\gamma_{\mathrm{max}}}}}^{+\infty} \frac{V_{1}-c^{-1}\gamma_{\mathrm{max}}}{V_{0}-c^{-1}\gamma_{\mathrm{max}}}\text{e}^{-cs}
        \dd s \notag\\
      &=& c^{-1}\left(1+\ln{\frac{V_{1}-c^{-1}\gamma_{\mathrm{max}}}{V_{0}-c^{-1}\gamma_{\mathrm{max}}}}\right)\notag
\end{eqnarray}
{which we term as $t_{d}$ in the context of this paper.}
Obviously, 
we have
\begin{equation}{\label{eq. new1}}
	\mathrm{P}\{\tau_{2i}-\tau_{2i-1}\geq s~|~x(\tau_{2i-1})\} \leq\mathrm{P}\left\{\hat{X}>s\right\},\quad \forall s\in \mathbb{R}
\end{equation}
{where  for $s\geq \frac{1}{c}\ln\frac{V_{1}-c^{-1}\gamma_{\mathrm{max}}}{V_{0}-c^{-1}\gamma_{\mathrm{max}}}$ cases  the inequality follows immediately from by \mbox{\eqref{Ineq distribution of down crossing time}}, while for $s< \frac{1}{c}\ln\frac{V_{1}-c^{-1}\gamma_{\mathrm{max}}}{V_{0}-c^{-1}\gamma_{\mathrm{max}}}$ cases the inequality follows from the upper boundedness of the distribution function.}

To cover the case of existing only finitely many loops in the noise-to-state exponentially stable system, we homoplastically construct a sequence of random variables $\left\{\hat{X}_{i}\right\}$ ($i\in\mathbb{Z}_{\textgreater 0}$) as
\begin{equation*}
   \hat{X_{i}}=
   \left\{
   \begin{array}{cl}
      \tau_{2i}-\tau_{2i-1}, & \tau_{2i-1}<+\infty; \\
      0,                     & \tau_{2i-1}=+\infty.
   \end{array}
   \right.
\end{equation*}
It is clear that $\mathrm{P}\left\{\hat{X}_{i}\geq s ~|~ \tau_{2i-1}=+\infty,\hat{X}_{1},\cdots,\hat{X}_{i-1}\right\}$ equals to $0$ if $s>0$, otherwise $1$, i.e., $\forall s\geq 0$ and $\forall i\in\mathbb{Z}_{\textgreater 0}$
\begin{equation}{\label{Ineq estimation of hat Xi}}
     \mathrm{P}\left\{\hat{X}_{i}\geq s ~|~ \tau_{2i-1}=+\infty,\hat{X}_{1},\cdots,\hat{X}_{i-1}\right\}
     \leq\mathrm{P}\left\{\hat{X}>s\right\}.
\end{equation}
{Moreover, for the case where $\tau_{2i-1}<\infty$, we have}
\begin{eqnarray}
	&&\mathrm{P}\left\{\hat{X}_{i}\geq s ~|~ \tau_{2i-1}<+\infty,\hat{X}_{1},\cdots,\hat{X}_{i-1}\right\} \notag\\
	&=&\text{E}\left[ \left. \mathrm{P}\left\{\hat{X}_{i}\geq s ~| ~x(\tau_{2i-1})\right\}~\right|~ \tau_{2i-1}<+\infty ,\hat{X}_{1},\cdots,\hat{X}_{i-1}\right] \notag \\
	&\leq&\mathrm{P}\left\{\hat{X}>s\right\}  \notag
\end{eqnarray}
{where the inequality follows immediately from \mbox{\eqref{eq. new1}}.}
This inequality together with the \text{inequality}~(\ref{Ineq estimation of hat Xi}) implies that that
\begin{equation}{\label{eq important 2}}
     \mathrm{P}\left\{\hat{X}_{i}\geq s ~|~\hat{X}_{1},\cdots,\hat{X}_{i-1}\right\}
     \leq \mathrm{P}\left\{\hat{X}>s\right\}, \quad \forall i\in\mathbb{Z}_{\textgreater 0} \text{ and } \forall  s\in \mathbb{R}.
\end{equation}
Hence, the series $\{\hat{X}_{i}\}$ ($i\in\mathbb{Z}_{\textgreater 0}$) is controlled by $\hat{X}$.

\begin{remark}
Proposition \mbox{\ref{Lem distribution of down crossing time}} and the above analysis imply that NSES can be connected with positive recurrence. Let $\tilde \tau_y$ be the first passage time when the trajectory with initial state $y$ hits $\{x|V(x)\leq V_0\}$. Then, we can similarly show that $\mathrm{P}\left({\tilde \tau_{y}}\geq s\right)\leq \frac{V(y)-c^{-1}\gamma_{\mathrm{max}}}{V_0-c^{-1}\gamma_{\mathrm{max}}} \text{e}^{-cs}$ for $V(y)>V_0$ and $\mathrm{P}\left({\tilde\tau_{y}}\geq s\right)=0$ otherwise. (The proof is almost the same as the one of Proposition \mbox{\ref{Lem distribution of down crossing time}}.)	Moreover, since the probability distribution of the down crossing time is controlled by an exponential function, the random variable $\tilde \tau_y$ is integrable, which implies an NSES system to be finite recurrent with respect to a compact set $\{x|V(x)\leq V_0\}$.
 \end{remark}

\section{Boundedness analysis of trajectory in the time-domain sense}{\label{Section Boundedness analysis}}
The above results indicate that the up crossing time and the down crossing time of every loop are both uniformly controlled by the corresponding random variables.
In this section, we will estimate the time bounds for loops to crossing up and down.
To be specific, the infimum/supremum of the average of up/down crossing time is estimated, based on which the lower bound of resident time ratio function may be reached.
\subsection{Infimum estimation of the average of the loop up crossing time} We take into key account that the number of loops is infinite. In this case, the average of the up crossing time can be expressed as $\lim_{n\to\infty}\sum_{i=1}^{n} \frac{\tilde{X}_{i}}{n}$. The form of this expression naturally motivates us to associate its estimation with the strong law of large numbers \cite{etemadi1981elementary} that says for a sequence of random variables $\{X_i\}$, $i\in\mathbb{Z}_{\textgreater 0}$, with independent identical distribution, if $\mathrm{E}\{X_1\}\textless +\infty$, then
   \begin{equation*}
      \lim_{n\to +\infty}\sum_{i=1}^{n} \frac{X_{i}}{n}=\text{E}\{X_1\}.
   \end{equation*}

To facilitate the applications of the strong law of large numbers, we manage to map as well as inversely map arbitrary random variables to independent identically distributed ones.

\begin{lemma}{\label{Lem random variables to uniform random variables}}
Assume $\{X_i\}_{i\in\mathbb{Z}_{\textgreater 0}}$ to be a sequence of random variables admitting the conditional distribution
\begin{eqnarray}
g\left(s~|~X_{1},\cdots,X_{i-1}\right)&=&
\mathrm{P}\left\{X_{i}\leq s~|~X_{1},\cdots,X_{i-1}\right\},\notag\\
g\left(s^{-}~|~X_{1},\cdots,X_{i-1}\right)
&=&\lim_{l\to s^-} g\left(l~|~X_{1},\cdots,X_{i-1}\right),\notag
\end{eqnarray}
where $s\in\mathbb{R}$, and $\{\xi_i\}_{i\in\mathbb{Z}_{\textgreater 0}}$ to be another sequence of random variables admitting the independent identical distribution $U(0,1)$, denoted by $\forall i\geq 1,~\xi_{i}\sim U(0,1)$. Furthermore, these two sequences are independent of each other. Then the series $\{Y_i\}_{i\in\mathbb{Z}_{\textgreater 0}}$ defined by
\begin{equation}{\label{eq. Y_{n}}}
    Y_{i}=\xi_{i}\cdot g\left(X_{i}~|~X_{1},\cdots,X_{i-1}\right)+(1-\xi_{i})\cdot g\left(X_{i}^{-}~|~X_{1},\cdots,X_{i-1}\right)
\end{equation}
is mutually independent and $\forall i\geq 1,~Y_{i}\sim U(0,1)$.
\end{lemma}

\begin{proof}
	The detail of the proof is shown in the appendix.
\end{proof}

Likewise, any random variable can be also generated from an independent identical distributed random variable.
\begin{lemma}{\label{Lem uniformly distributed random variables to any random variables 1}}
Given a random variable $Y\sim U(0,1)$, let $F(\cdot)$ be any distribution function, the random variable $Z$ is subject to the distribution $F(\cdot)$ if it takes any one of the following definitions with $s\in\mathbb{R}$
\begin{enumerate}
  \item[\rm{(\romannumeral1)}] $Z:=\inf\{s~|~F(s)\geq Y\}$;
  \item[\rm{(\romannumeral2)}] $Z:=\sup\{s~|~F(s) \leq Y\}$.
\end{enumerate}
\end{lemma}

\begin{proof}
	(i) In the case of $Z:=\inf\{s~|~F(s)\geq Y\}$, since $F(\cdot)$ is right continuous, there is $F(Z)=\lim_{s\to Z^{+}}F(s)\geq Y$. This suggests that $Z$ is the minimum point of the set $\{s~|~F(s)\geq Y\}$. Note that $F(\cdot)$ has monotonicity, so we have $\{s \geq Z\}=\{F(s)\geq Y\}$. Therefore,
	\begin{equation*}
	\mathrm{P}\{Z\leq s\}= \mathrm{P}\{  Y\leq F(s)\}=F(s),
	\end{equation*}
	which indicates that $Z:=\inf\{s~|~F(s)\geq Y\}$ is subject to the distribution $F(\cdot)$.
	
	(ii) In the case of $Z:=\sup\{s~|~F(s) \leq Y\}$, we have $F\left(Z^{-}\right)\leq Y$. The monotonicity of $F(\cdot)$ also suggests $\{Z \geq s\}=\{F(s^{-})\leq Y\}$. We then get
	\begin{equation*}
	\mathrm{P}\{Z < s\}= 1 - \mathrm{P}\{ Z \geq s \}= 1-\mathrm{P}\{Y>F(s^{-})\}=F(s^{-}).
	\end{equation*}
	This completes the proof.
\end{proof}

Utilizing above mapping methods, we can estimate the infimum average of a sequence of {random variables}, which is not necessarily identically distributed {and also not necessarily independent}.

\begin{theorem}{\label{thm SLLN lower bound esitmation}}
	Let $\left\{X_{i}\right\}_{i\in\mathbb{Z}_{\textgreater 0}}$ be a sequence of random variables and $X$ be another random variable satisfying $\mathrm{E}\{ X\}< +\infty$.
	If for all $i\in\mathbb{Z}_{\textgreater 0}$ and $s\in\mathbb{R}$ there is
	\begin{equation}{\label{Ineq heavy tail condition}}
	\mathrm{P}\left\{X_{i}>s~|~X_{1},\cdots,X_{i-1}\right\} \geq \mathrm{P} \{X>s\},
	\end{equation}
	then
	\begin{equation}\label{Ineq heavy tail condition1}
	\liminf_{n\to+\infty}\sum_{i=1}^{n} \frac{X_{i}}{n} \geq \mathrm{E}\{ X\} \quad \quad \mathrm{a.~s.}
	\end{equation}
\end{theorem}

\begin{proof}
Let $\{\xi_{i}\}_{i\in\mathbb{Z}_{\textgreater 0}}$ be a sequence of random variables satisfying $\xi_{i}\sim U(0,1)$, and moreover, $\{\xi_{i}\}$ is independent of $\{X_{i}\}$. In the meanwhile, define $Y_{i}$ as Eq. (\ref{eq. Y_{n}}) and denote
	\begin{equation}{\label{eq. Z2}}
	Z_{i}=\inf\{s~|~F(s)\geq Y_{i}\}
	\end{equation}
where $F(s)=\mathrm{P}\{X\leq s\}$. According to \textit{Lemma \ref{Lem random variables to uniform random variables}}, $\{Y_{i}\}_{i\in\mathbb{Z}_{\textgreater 0}}$ is mutually independent and $Y_{i}\sim U(0,1)$.
	Thus, $\{Z_{i}\}_{i\in\mathbb{Z}_{\textgreater 0}}$ is also mutually independent and, by \textit{Lemma \ref{Lem uniformly distributed random variables to any random variables 1}}, $Z_{i}$ is subject to $F(\cdot)$. Hence, based on the strong law of large numbers \cite{etemadi1981elementary}, we have
	\begin{equation*}
	\lim_{n\to\infty}\sum_{i=1}^{n} \frac{Z_{i}}{n}
	=\mathrm{E} \{Z_{1}\} =\text{E} \{X\} \quad \quad \mathrm{a.~s.}
	\end{equation*}
	
Further, from the inequality (\ref{Ineq heavy tail condition}), we get $g(s|X_{1},\cdots,X_{i-1})\leq F(s)$. Combing this result with the fact $Y_{i}\leq g(X_{i}~|~X_{1},\cdots,X_{i-1})$ (cf. Eq. (\ref{eq. Y_{n}})) yields
	\begin{equation*}
	F(X_{i}^{}) \geq g(X_{i}~|~X_{1},\cdots,X_{i-1})\geq Y_{i}.
	\end{equation*}
	Thus, by {above relation and} Eq. (\ref{eq. Z2}), we know for all $i\in\mathbb{Z}_{\textgreater 0}$. So there is
	\begin{equation*}
	\liminf_{n\to\infty} \sum_{i=1}^{n} \frac{X_{i}}{n}
	\geq \liminf_{n\to\infty} \sum_{i=1}^{n} \frac{Z_{i}}{n}=\text{E}\{ X\} \quad\quad \mathrm{a.~s.}
	\end{equation*}
The desired result is obtained.
\end{proof}

Applying \textit{Theorem} \ref{thm SLLN lower bound esitmation} to the loop up crossing time series $\{\tau_{2i+1}-\tau_{2i}\}_{i\in\mathbb{Z}_{\geq 0}}$, we get the estimation of the infimum of the average of the series.

\begin{proposition}{\label{Lem up crossing time is almost surely lower bounded}}
For a noise-to-state exponentially stable system in the form of Eq. (\ref{Eq. stochastic systems}), let $\{\tau_{2i+1}-\tau_{2i}\}_{i\in\mathbb{Z}_{> 0}}$ represent the loop up crossing time series and $t_{\mathrm{u}}$ be the expectation of the random variable $\tilde{X}$ that controls uniformly the series, then we have
	\begin{equation*}
	 \mathrm{P}\left\{
	 \left\{\text{\textup{There are \text{infinitely many loops}}}\right\}
	 \bigcap\left\{
	 \liminf_{n\to+\infty} \sum_{i=1}^{n} \frac{\tau_{2i+1}-\tau_{2i}}{n} {\textless} t_{\mathrm{u}}
	 \right\}
	 \right\}=0.
	\end{equation*}
\end{proposition}

\begin{proof}
The result is straightforward since
	\begin{eqnarray}
		 &&\mathrm{P}\left\{
		 \left\{\text{There are \text{infinitely many loops}}\right\}
		 \bigcap\left\{
		 \liminf_{n\to+\infty} \sum_{i=1}^{n} \frac{\tau_{2i+1}-\tau_{2i}}{n} {\textless} t_{\mathrm{u}}
		 \right\}
		 \right\}\notag\\
		 &=&\mathrm{P}\left\{
		 \left\{\text{There are \text{infinitely many loops}}\right\}
		 \bigcap\left\{
		 \liminf_{n\to+\infty} \sum_{i=1}^{n} \frac{\tilde{X}_{i}}{n} {\textless} t_{\mathrm{u}}
		 \right\}
		 \right\}\notag\\
		 &\leq& \mathrm{P}\left\{
		 \liminf_{n\to+\infty} \sum_{i=1}^{n} \frac{\tilde{X}_{i}}{n} {\textless} t_{\mathrm{u}}
		 \right\}\notag
		 =0,
	\end{eqnarray}
where the last equality holds due to Eq. \eqref{eq. important 1} and (\ref{Ineq heavy tail condition1}).
\end{proof}

The above result implies that the inferior average of the loop up crossing time is almost surely lower bounded by $t_\mathrm{u}$ if the number of the loops is infinite.

\begin{remark}{\label{remark 4}}
	{Also, this result indicates $\lim_{i\to\infty} \tau_{i}$, a term that is greater than or equals to the infinite sum of $\tilde{X}_{i}$, to be almost surely infinite. Therefore, the decomposition of the whole evolution time through each loop, as mentioned in \textit{Remark 2}, is well defined.  }
\end{remark}

\subsection{Supremum estimation of the average of the loop down crossing time}
The supremum of the average of the loop down crossing time can be also estimated based on the same principle as applied to the loop up crossing time.

\begin{theorem}{\label{thm SLLN upper bound estimation}}
Assume that $\left\{X_{i}\right\}_{i\in\mathbb{Z}_{\textgreater 0}}$ is a sequence of random variables and $X$ is another random variable satisfying $\mathrm{E}\{X\}<+\infty$. If $\forall i\in\mathbb{Z}_{\textgreater 0}$ and $\forall s\in\mathbb{R}$ there is
	\begin{equation}{\label{Ineq light tail condition}}
	\mathrm{P}\left\{X_{i}\geq s~|~X_{1},\cdots,X_{i-1}\right\}\leq \mathrm{P} \{X>s\}
	\end{equation}
	then
	\begin{equation}\label{Ineq light tail condition1}
	\limsup_{n\to\infty}\sum_{i=1}^{n} \frac{X_{i}}{n} \leq \mathrm{E} \{X\} \quad \quad \mathrm{a.~s.}
	\end{equation}
\end{theorem}

\begin{proof}
Similar to the proof of \textit{Theorem} \ref{thm SLLN lower bound esitmation}, $\forall i\in\mathbb{Z}_{\textgreater 0}$ we construct $Y_{i}$ according to Eq. (\ref{eq. Y_{n}}) and define $Z_{i}$ by
	\begin{equation}{\label{eq. Z1}}
	 Z_{i}=\mathrm{sup}\{s~|~F(s)\leq Y_{i}\}
	\end{equation}
where $F(s)=\mathrm{P}\{X\leq s\}$. Naturally, $\{Z_{i}\}_{i\in\mathbb{Z}_{\textgreater 0}}$ is mutually independent and every $Z_{i}$ is subject to the distribution $F(\cdot)$. Hence, we also have
	\begin{equation*}
	  \lim_{n\to\infty}\sum_{i=1}^{n} \frac{Z_{i}}{n}
	  =\text{E} \{Z_{1}\} =\text{E} \{X\} \quad \quad \mathrm{a.~s.}
	\end{equation*}
	
Besides, the inequality (\ref{Ineq light tail condition}) implies $g(s^{-}|X_{1},\cdots,X_{i-1})\geq F(s)$.
Therefore, we get
	\begin{equation*}
	 F(X_{i}^{}) \leq g(X_{i}^{-}~|~X_{1},\cdots,X_{i-1}) \leq Y_{i}.
	\end{equation*}
This combining Eq. (\ref{eq. Z1}) means $X_{i}\leq Z_{i}$ for all $i\in\mathbb{Z}_{\textgreater 0}$. So the result expressed as
	\begin{equation*}
	 \limsup_{n\to\infty} \sum_{i=1}^{n} \frac{X_{i}}{n}
	 \leq \limsup_{n\to\infty} \sum_{i=1}^{n} \frac{Z_{i}}{n}=\text{E}\{ X\} \quad\quad \mathrm{a.~s.}
	\end{equation*}
is true.
\end{proof}

We apply this theorem to the loop down crossing time series and can obtain the supremum estimation of their average value.

\begin{proposition}{\label{Lem down crossing time is almost surely upper bounded}}
Assume that a stochastic system governed by (\ref{Eq. stochastic systems}) has NSES, then we have
	\begin{equation*}
	\mathrm{P}\left\{
	\left\{\text{\textup{There are \text{infinitely many loops}}}\right\}
	\bigcap\left\{
	\limsup_{n\to\infty} \sum_{i=1}^{n} \frac{\tau_{2i}-\tau_{2i-1}}{n} {\textgreater} t_{\mathrm{d}}
	\right\}
	\right\}=0,
	\end{equation*}
where $\{\tau_{2i}-\tau_{2i-1}\}_{i\in\mathbb{Z}_{\textgreater 0}}$ is the loop down crossing time series, and $t_{\mathrm{d}}$ {is the expectation of the random variable $\hat{X}$ that controls the series uniformly}.
\end{proposition}

\begin{proof}
	It is simply because
	\begin{eqnarray}
	&&\mathrm{P}\left\{
	\left\{\text{There are \text{infinitely many loops}}\right\}
	\bigcap\left\{
	\limsup_{n\to+\infty} \sum_{i=1}^{n} \frac{\tau_{2i}-\tau_{2i-1}}{n} {\textgreater} t_{\mathrm{d}}
	\right\}
	\right\}\notag\\
	&=&\mathrm{P}\left\{
	\left\{\text{There are \text{infinitely many loops}}\right\}
	\bigcap\left\{
	\limsup_{n\to+\infty} \sum_{i=1}^{n} \frac{\hat{X}_{i}}{n} {\textgreater} t_{\mathrm{d}}
	\right\}
	\right\}\notag\\
	&\leq& \mathrm{P}\left\{
	\limsup_{n\to+\infty} \sum_{i=1}^{n} \frac{\hat{X}_{i}}{n} {\textgreater} t_{\mathrm{d}}
	\right\}\notag
	=0.
	\end{eqnarray}
The last equality applies to Eq. \eqref{eq important 2} and (\ref{Ineq light tail condition1}).
\end{proof}

We can conclude from this proposition that the superior average of the loop down crossing time is almost surely upper bounded by $t_\mathrm{d}$ if the number of the loops is infinite.

\subsection{Lower bound estimation of trajectory all through evolution} The estimations of the loop up and down crossing time bounds for the noise-to-state exponentially stable system essentially reflect the ratio of resident time. Based on these estimations, it is possible to derive time lower bound for the trajectories to evolve with a relatively ``low energy region". The following theorem gives this time lower bound.

\begin{theorem}{\label{thm time average analysis of V}}
	For a noise-to-state exponentially stable system (\ref{Eq. stochastic systems}), we have
	\begin{equation*}
	\liminf_{T\to+\infty} \frac{1}{T}\int_{0}^{T} \mathbbold{1}_{\left\{V(x(s)) < V_{1}\right\}} d s
	\geq \frac{t_{\mathrm{u}}}{t_{\mathrm{u}}+t_{\mathrm{d}}}\quad\quad {\mathrm{a.~s.}},
	\end{equation*}
where $V(x)$ is the noise-to-state Lyapunov function, $V_1$ shares the same meaning as in \textit{Definition} \ref{Def loops}, and $t_{\mathrm{u}}$ and $t_{\mathrm{d}}$ are the same as in \textit{Proposition} \ref{Lem up crossing time is almost surely lower bounded} and \ref{Lem down crossing time is almost surely upper bounded}, respectively.
\end{theorem}

\begin{proof}
	The detail of the proof is shown in the appendix.
\end{proof}

As a lower bound, it is natural to expect this value, i.e., ${t_{\mathrm{u}}}/(t_{\mathrm{u}}+t_{\mathrm{d}})$, as big as possible.
According to \textit{Definition} \ref{Def loops}, the loop up and down crossing time closely depends on the choice of $V_{0}$ when $V_1$.
Hence, here comes an idea to optimize the parameter $V_0$ so that $\frac{t_{\mathrm{u}}}{t_{\mathrm{u}}+t_{\mathrm{d}}}$ reaches the maximum or equivalently $\frac{t_{\mathrm{d}}}{t_{\mathrm{u}}}$ reaches the minimum.

\begin{lemma}\label{LemGao}
	For a noise-to-state exponentially stable system (\ref{Eq. stochastic systems}), under a given $V_1$ there is
	\begin{equation*}
	\min_{V_0} \frac{t_{\mathrm{d}}}{t_{\mathrm{u}}}=\frac{-\mathcal{W}_{-1}(-\mathrm{e}^{-2})}{ \ln{\frac{V_{1}}{c^{-1}\gamma_{\mathrm{max}}}}}
	\end{equation*}
	where $V_0$ shares the same meaning as the one in \textit{Definition} \ref{Def loops}, and $V_1$, $t_{\mathrm{u}}$ and $t_{\mathrm{d}}$ are the same as the ones in \textit{Theorem} \ref{thm time average analysis of V}.
\end{lemma}
\begin{proof}
	Denote $\beta=\frac{V_{0}-c^{-1}\gamma_{\mathrm{max}}}{V_{1}-c^{-1}\gamma_{\mathrm{max}}}$, then under a given $V_1$ there is $\min_{V_0} \frac{t_{\mathrm{d}}}{t_{\mathrm{u}}}=\min_{\beta} \frac{t_{\mathrm{d}}}{t_{\mathrm{u}}}$. From the expressions of $t_{\mathrm{u}}$ and $t_{\mathrm{d}}$, i.e., Eqs. (\ref{eq calculate t_{uc}}) and (\ref{eq calculate t_{dc}}), we have
	\begin{equation*}
	t_{\mathrm{u}} = (1-\beta)c^{-1}\ln{\frac{V_{1}}{c^{-1}\gamma_{\mathrm{max}}}} \quad \text{and} \quad  t_{\mathrm{d}}= c^{-1}(1-\ln{\beta}).
	\end{equation*}
	Further, we get
	\begin{equation*}
	\frac{\dd~ \frac{t_{\mathrm{d}}}{t_{\mathrm{u}}}}{d~ \beta}
	=\frac{2-\frac{1}{\beta}-\ln{\beta}}
	{(1-\beta)^{2}\ln{\frac{V_{1}}{c^{-1}\gamma_{\mathrm{max}}}}}
	=\frac{\ln\left(\text{e}^{2}\cdot\frac{1}{\beta}\text{e}^{-\frac{1}{\beta}}\right)}            {(1-\beta)^{2}\ln{\frac{V_{1}}{c^{-1}\gamma_{\mathrm{max}}}}}.
	\end{equation*}
	Note that $\beta\in(0,1)$.
	When $0<\beta<\frac{-1}{\mathcal{W}_{-1}(-\text{e}^{-2})}$ there is $\frac{\dd~ \frac{t_{\mathrm{d}}}{t_{\mathrm{u}}}}{\dd~ \beta}\textless 0$, while when $\frac{-1}{\mathcal{W}_{-1}(-\text{e}^{-2})}<\beta<1$ there is $\frac{\dd~ \frac{t_{\mathrm{d}}}{t_{\mathrm{u}}}}{\dd~ \beta}\textgreater 0$. This implies that $\frac{t_{\mathrm{d}}}{t_{\mathrm{u}}}$ reaches the minimum $\frac{-\mathcal{W}_{-1}(-\text{e}^{-2})}{ \ln{\frac{V_{1}}{c^{-1}\gamma_{\mathrm{max}}}}}$ at $\beta=\frac{-1}{\mathcal{W}_{-1}(-\text{e}^{-2})}$. It should be pointed out that \textit{Definition} \ref{Def loops} (loop) requests $c^{-1}\gamma_{\mathrm{max}}\textless V_0\textless V_1$. According to the optimal $\beta=\frac{-1}{\mathcal{W}_{-1}(-\text{e}^{-2})}$, the corresponding $V_0$ equals to $\frac{V_1-c^{-1}\gamma_{\mathrm{max}}}{-\mathcal{W}_{-1}(-\text{e}^{-2})}+c^{-1}\gamma_{\mathrm{max}}$, which obviously satisfies the condition $c^{-1}\gamma_{\mathrm{max}}\textless V_0\textless V_1$.
\end{proof}

Following the above theorem and lemma we can provide a lower bound for the resident time ration function $D(r)$.

\begin{theorem}{\label{thm main result}}
	For a noise-to-state exponentially stable system governed by Eq. (\ref{Eq. stochastic systems}), for any $r\textgreater \alpha^{-1}_{1}(c^{-1}\gamma_{\mathrm{max}})$ we have
	\begin{eqnarray}{\label{Ineq. main result}}
	D(r)\triangleq
	\liminf_{T\to\infty} \frac{1}{T}\int_{0}^{T} \mathbbold{1}_{\left\{\left\|x(s)\right\|_{2} < r \right\}} \dd s
	&\geq& 		
	\frac{\ln{\frac{\alpha_{1}(r)}{c^{-1}\gamma_{\mathrm{max}}}}}
	{ -\mathcal{W}_{-1}(-\text{e}^{-2}) + \ln{\frac{\alpha_{1}(r)}{c^{-1}\gamma_{\mathrm{max}}}}} \quad\quad {\mathrm{a.~s.}}
	\end{eqnarray}
\end{theorem}

\begin{proof}
	According to \textit{Definition} \ref{Def NNS}, there is $\mathbbold{1}_{\left\{\left\|x(s)\right\|_{2}<r\right\}}=\mathbbold{1}_{\left\{\alpha_1(\left\|x(s)\right\|_{2}) < \alpha_1(r)\right\}}$. Since $\alpha_1(\left\|x(s)\right\|_{2})\leq V(x(s))$, we get
	\begin{eqnarray}\label{gao}
	D(r)		&\geq&
	\liminf_{T\to\infty} \frac{1}{T}\int_{0}^{T} \mathbbold{1}_{\left\{V(x(s)) < \alpha_1(r)\right\}} \dd s.
	\end{eqnarray}
	From \textit{Lemma} \ref{LemGao}, with the minimal $\frac{t_{\mathrm{d}}}{t_{\mathrm{u}}}$ there is
	\begin{equation*}
	\liminf_{T\to\infty} \frac{1}{T}\int_{0}^{T} \mathbbold{1}_{\left\{V(x(s)) < V_{1}\right\}} \dd s
	\geq
	\frac{\ln{\frac{V_{1}}{c^{-1}\gamma_{\mathrm{max}}}}}
	{ -\mathcal{W}_{-1}(-\text{e}^{-2}) + \ln{\frac{V_{1}}{c^{-1}\gamma_{\mathrm{max}}}}}\quad\quad {\mathrm{a.~s.}}.
	\end{equation*}
	Since $\forall r\textgreater \alpha^{-1}_{1}(c^{-1}\gamma_{\mathrm{max}})$, there is $\alpha_1(r)\textgreater c^{-1}\gamma_{\mathrm{max}}$. Replacing $V_1$ by $\alpha_1(r)$ in the above inequality and further combining the inequality (\ref{gao}), we get the inequality (\ref{Ineq. main result}).
\end{proof}

\begin{remark}{\label{Rem b}}
	Since $\alpha_{1}(\cdot)\in\mathcal{K}_{\infty}$ as shown in \textit{Definition \ref{Def NNS}}, the lower bound function (with respect to $r$) of $D(r)$ in the inequality (\ref{Ineq. main result}) satisfies
	$$\lim_{r\to+\infty}\frac{\ln{\frac{\alpha_{1}(r)}{c^{-1}\gamma_{\mathrm{max}}}}}
	{ -\mathcal{W}_{-1}(-\text{e}^{-2}) + \ln{\frac{\alpha_{1}(r)}{c^{-1}\gamma_{\mathrm{max}}}}}=1.$$
	This means that the above lower bound function is a good candidate for $b(r)$ in Eq. (\ref{Eq. goal0}), i.e.,
	\begin{equation}{\label{Eq. b}}
	b(r)=\frac{\ln{\frac{\alpha_{1}(r)}{c^{-1}\gamma_{\mathrm{max}}}}}
	{ -\mathcal{W}_{-1}(-\text{e}^{-2}) + \ln{\frac{\alpha_{1}(r)}{c^{-1}\gamma_{\mathrm{max}}}}}, \quad\quad \forall r\textgreater \alpha^{-1}_{1}(c^{-1}\gamma_{\mathrm{max}}).
	\end{equation}
Denote $q_{k}(\gamma_{\mathrm{max}})$ as the $k$ fractile of $b(r)$ with respect to the maximum noise intensity $\gamma_{\mathrm{max}}$, then according to $b(q_k(\gamma_{\mathrm{max}}))=k$ we can calculate $q_{k}(\gamma_{\mathrm{max}})$ to be
    \begin{equation}{\label{Eq. qk}}
       q_{k}(\gamma_{\mathrm{max}})=\alpha^{-1}_{1}
            \left(         c^{-1}\gamma_{\mathrm{max}}\cdot\exp\left(-\frac{k}{1-k}\mathcal{W}_{-1}\left(-\text{e}^{-2}\right)\right)
            \right),
       \quad\quad k\in(0,1).
    \end{equation}
which naturally satisfies $q_k(\gamma_{\mathrm{max}})\textgreater \alpha^{-1}_{1}(c^{-1}\gamma_{\mathrm{max}})$. Based on Eq. (\ref{Eq. qk}), it can be concluded that for any given $k\in(0,1)$ and any trajectory of a noise-to-state exponentially stable system, there is more than $100\cdot k$ percent of all time for the trajectory to evolve in a closed ball $\bar{\mathbb{B}}^n(q_k)$.
\end{remark}

Therefore, similar to boundedness in probability, we can present ``boundedness in the time horizon" as follow.

\begin{theorem}{\label{thm boundedness in time horizon}}
	For a noise-to-state exponentially stable system governed by Eq. (\ref{Eq. stochastic systems}),
	for any $\epsilon>0$ there is a $\mathcal{K_{\infty}}$ function $\delta$ such that
	\begin{equation}\label{eq. boundedness in the time horizon}
	D\left(\delta\left(\sup_{t}\|\Sigma(t)\|_{F}\right)\right)
	\geq 1-\epsilon
    \quad\quad {\mathrm{a.~s.}}
	\end{equation}
	where $\delta(\cdot)=q_{1-\epsilon}(\cdot)$.
\end{theorem}

\begin{proof}
	The equality holds immediately by \eqref{Ineq. main result} \eqref{Eq. b} and \eqref{Eq. qk}, and the $\mathcal{K_{\infty}}$ property of $q_{1-\epsilon}(\cdot)$ follows immediately from $\alpha_{1}$ being $\mathcal{K_{\infty}}$.
\end{proof}

\begin{remark}
	Since, for a fixed trajectory, $D(r)$ is nondecreasing, right continuous, and upper bounded by one, we can also call $D(r)$ a distribution-like function.
    Therefore, we can observe that \eqref{eq. boundedness in the time horizon} and \eqref{eq. probability nss} share the same structure.
    Moreover, similar to \eqref{eq. probability nss} which provides a region to find the state with a given probability,  inequality \eqref{eq. boundedness in the time horizon} provides another region where the trajectory evolutes with a given proportion of time.
    In this regard, the ``boundedness in the time horizon", \eqref{eq. boundedness in the time horizon} is a good extension of the boundedness in probability \eqref{eq. probability nss} from a time domain view.
\end{remark}

\begin{remark}
	The lower bound function $b(r)$ is influenced by $\gamma_{\mathrm{max}}$, where the latter one is evaluated by the noise intensity $\Sigma(t)$.
	When the noise intensity decays, i.e., $\gamma_{\mathrm{max}}$ decreases, the lower bound estimation, $b(r)$, will become larger. Particularly, when $\Sigma(t)=\boldmath{0}$, the resident time ratio function strictly equals to one, which means trajectories will stay in a bounded region forever.
	Conversely, when the noise intensity grows stronger, $b(r)$ decreases and therefore the ``stability" will be degenerated.
	This indicates that it is difficult for a ``strong" disturbed system to be stable in a bounded region for a long time.
	Such conclusion can also be obtained from $q_{k}$ and \textit{Theorem \ref{thm boundedness in time horizon}}.
	These results provide an additional insight on understanding the concept of the noise-to-state stability from a point of time-domain view.
\end{remark}
\section{A numerical example}{\label{Section numerical example}}
\begin{figure}
	\centering
	\includegraphics[width=1\textwidth]{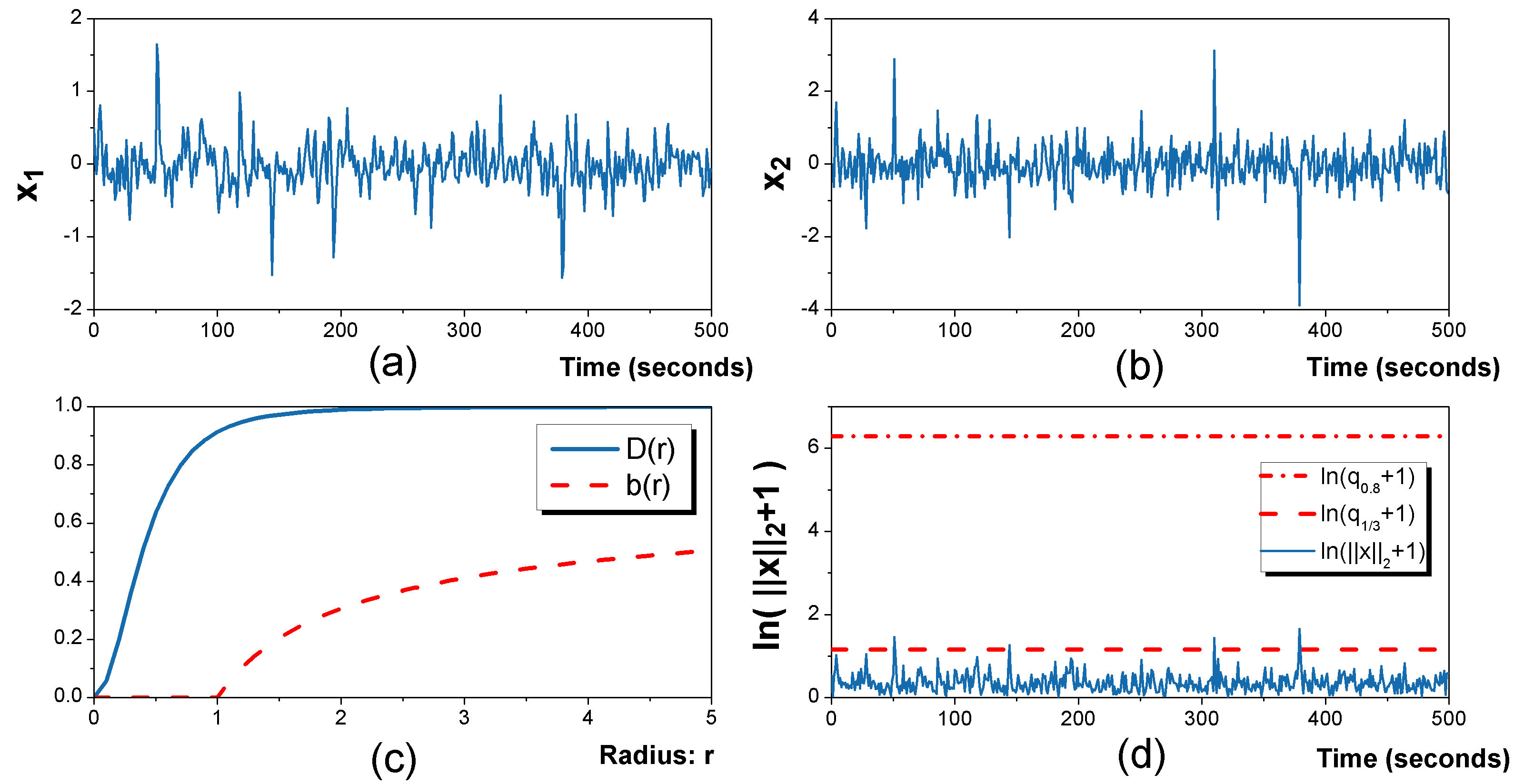}\\
	\caption{The evolution behaviors of the system (\ref{Eq numerical example system}): (a) the evolution of $x_1(t)$; (b) the evolution of $x_2(t)$; (c) the distribution function $D(r)$ versus its lower bound function $b(r)$;
	(d) the evolution of the state norm versus $q_{1/3}$ and $q_{0.8}$ drawn in the logarithmic form.}
	\label{Fig simulation}
\end{figure}
In this section, we present a numerical example to illustrate our result.

Consider a system in the form of
\begin{equation}{\label{Eq numerical example system}}
\dd \left(
\begin{array}{c}
x_{1}\\
x_{2}
\end{array}
\right)=
\underbrace{\left(
	\begin{array}{c}
	-x_{1}+x_{2}\\
	-x_{1}-x_{2}
	\end{array}
	\right)}_{f(x)} \dd t
+
\underbrace{\left(
	\begin{array}{cc}
	0 & 0\\
	x_{2} & 1
	\end{array}
	\right)}_{h(x)}
\underbrace{\left(
	\begin{array}{cc}
	1 & 0\\
	0 & \sin{t}
	\end{array}
	\right)}_{\Sigma(t)}
\dd \left(
\begin{array}{c}
\mathcal{B}_{1}\\
\mathcal{B}_{2}
\end{array}
\right).
\end{equation}
Obviously, its structure consists with the system (\ref{Eq. stochastic systems}). Taking
a positive definite function $V(x)=\frac{1}{2}\left(x_{1}^{2}+x_{2}^{2}\right)$ and substituting it into Eq. (\ref{Eq. def NSS}), we have
\begin{equation*}
\mathcal{L}\left[V(x)\right]
=-\left(x_{1}^{2}+x_{2}^{2}\right)+\frac{1}{2}\left(x_{2}^{2}+\sin^{2}{t}\right)
\leq -\frac{1}{2}\left(x_{1}^{2}+x_{2}^{2}\right)
+\frac{1}{2}\sqrt{\left\|\Sigma(t)\right\|^{2}_{F}-1}.
\end{equation*}
If we define $\alpha_{1}\left(\|x\|_{2}\right)=\alpha_{2}\left(\|x\|_{2}\right)=\alpha_{3}\left(\|x\|_{2}\right)=\frac{1}{2}\|x\|^{2}_{2}$, $c=1$ and $\gamma\left(\|\Sigma(t)\|_{F}\right)=\frac{1}{2}\sqrt{\left\|\Sigma(t)\right\|^{2}_{F}-1}$, then combing the above inequality we can conclude
\begin{equation*}
\alpha_{1}\left(\|x\|_{2}\right) = V(x) = \alpha_{2}\left(\|x\|_{2}\right)
\end{equation*}
and
\begin{equation*}
\mathcal{L}\left[V(x)\right]\leq
-\alpha_{3}\left(\|x\|_{2}\right)+\gamma\left(\|\Sigma(t)\|_{F}\right)
= -cV(x)+\gamma\left(\|\Sigma(t)\|_{F}\right),
\end{equation*}
which are in accordance with the conditions (\ref{Eq. control of V}) and (\ref{Eq. def NSS}) in \textit{Definition \ref{Def NNS}}.
Therefore, the system (\ref{Eq numerical example system}) is a noise-to-state exponentially stable system with the noise-to-stable Lyapunov function to be $V(x)=\frac{1}{2}\left(x_{1}^{2}+x_{2}^{2}\right)$.
Moreover, $\gamma_{\max}=\frac{1}{2}$ for this system. Further, we can verify that the noise term is singular, which will lead to the classic ergodicity analysis invalid for this system. However, it can be made time-domain analysis with respect to NSES using our method. Apply \textit{Theorem \ref{thm main result}} to this example, then we get that for any $r>\alpha_{1}^{-1}\left(c^{-1}\gamma_{\mathrm{max}}\right)=1$
\begin{equation}\label{gao2}
D(r)\triangleq
\liminf_{T\to\infty} \frac{1}{T}\int_{0}^{T} \mathbbold{1}_{\left\{\left\|x(s)\right\|_{2} <  r\right\}} \dd s
\geq b(r) \triangleq		
\frac{2\ln{r}}
{ -\mathcal{W}_{-1}(-\text{e}^{-2}) + 2\ln{r}}\quad\quad {\mathrm{a.~s.}}
\end{equation}
Besides, according to \eqref{Eq. qk}, we obtain $q_{\frac{1}{3}}=\exp\left\{-0.25\mathcal{W}_{-1}\left(-e^{-2}\right)\right\}\approx 2.2$
and
$q_{0.8}=\exp\left\{-2\mathcal{W}_{-1}\left(-e^{-2}\right)\right\}\approx 540$.
Therefore, the state of the system \eqref{Eq numerical example system} stays in the closed ball $\bar{\mathbb{B}}^2(q_{{1}/{3}})$ (approximating $\bar{\mathbb{B}}^2(2.2)$) and $\bar{\mathbb{B}}^2(q_{0.8})$ (approximating $\bar{\mathbb{B}}^2(540)$) for more than one third of the time and $80$\% of the time, respectively.
In order to exhibit the evolution behaviors of the trajectory, Fig. \ref{Fig simulation} (a) and (b) report the time recordings of $x_1(t)$ and $x_2(t)$ within $500$ s, respectively. It is obvious that there are strong oscillating phenomena in the evolution process. The state basically fluctuates in an irregular way around the origin point $(0,0)^\top$, the equilibrium for the corresponding deterministic system.
According to the numerical experiment, the resident time ratio function, $D(r)$, is firmly controlled by $b(r)$ defined in Eq. (\ref{gao2}), which can be observed from Fig. \ref{Fig simulation} (c).
When $r\textgreater 1$, $D(r)$ is always larger than $b(r)$.
Furthermore, we give the evolution of $\|x(t)\|_2$ in the logarithmic form in Fig. \ref{Fig simulation} (d), from which it can be seen that the state is bounded by the closed ball $\bar{\mathbb{B}}^2(q_{{1}/{3}})$ and $\bar{\mathbb{B}}^2(q_{0.8})$ for more than (actually far more than) one third of the time and 80\% of the time, respectively. Once the state exceeds the ball boundary, it will soon be pulled back. All these observations tally with the preceding theoretical analysis well.
\section{Conclusion}{\label{Section conclusion}}
In this paper, we analyze noise-to-state exponentially stable systems from the time-domain perspective.
Mainly, we show the resident time ratio function to be almost surely lower bounded and present a particular kind of stability for noise-to-state exponentially stable systems. First, we construct loops for the trajectory and show the probability distribution of up/down crossing time of these loops is uniformly controlled.
Based on it, the average of up/down crossing time is estimated and the main result, boundedness in the time horizon, is presented.
In the meanwhile, we point out that as the noise intensity tends to 0 the resident time ratio function goes to 1 and therefore the system becomes more stable.
This result provides an interpretation of NSES from the time-domain perspective.

There are some future works that we are concerned with related to the time-domain interpretation of NSES. First, according to the numerical example, the lower bound function $b(\cdot)$ is not very accurate to approximate the resident time ratio function $D(\cdot)$, which call for more efforts of ours to improve it.
Besides, we also plan to extend these results to noise-to-state exponentially stable systems with other forms of noise, such as jump disturbances and switch behaviors.
\section*{Appendix}

In the appendix, we show some proofs of our results.

\begin{proof}[Proof of Lemma 3]
	As a distribution function, $g(\cdot)$ is non-decreasing and right continuous, which together with $\forall i\geq 1,~\xi_{i}\sim U(0,1)$ yields $$g\left(X_{i}~|~X_{1},\cdots,X_{i-1}\right)\geq Y_i\geq g\left(X_{i}^{-}~|~X_{1},\cdots,X_{i-1}\right).$$ Denote $L_i=\mathrm{sup}\{l~|~g\left(l~|~X_{1},\cdots,X_{i-1}\right) < s\}$, then we have
	\begin{equation*}
	g\left(L_i~|~X_{1},\cdots,X_{i-1}\right)
	\geq s
	\quad \text{and} \quad
	g\left(L_i^{-}~|~X_{1},\cdots,X_{i-1}\right)
	\leq s.
	\end{equation*}
	{We continue the proof in separate two cases.}
	
	{Cases I): $g\left(L_i~|~X_{1},\cdots,X_{i-1}\right)\geq s$ and $g\left(L_i^{-}~|~X_{1},\cdots,X_{i-1}\right)\textless s$ for a given $s\in(0,1)$, i.e., the former larger than the latter. For simplicity of notations, denote $X_{1},\cdots,X_{i-1}$ and $\xi_{1},\cdots,\xi_{i-1}$ by $X_{1:~i-1}$ and $\xi_{1:~i-1}$, respectively. Then with $s\in(0,1)$ we consider the conditional probability $\mathrm{P}\left\{Y_{i} < s~|~X_{1:~i-1},\xi_{1:~i-1}\right\}$ that follows}
	\begin{eqnarray}\label{Lemma3gao}
	&~&\mathrm{P}\left\{Y_{i} < s~|~X_{1:~i-1},\xi_{1:~i-1}\right\}\\
	&=& {\mathrm{P}\left\{
		g\left(X_{i}~|~X_{1:~i-1}\right) < s
		~|~X_{1:~i-1}\right\}} \notag\\
	&&+{\mathrm{P}\left\{
		g\left(X_{i}~|~X_{1:~i-1}\right)\geq s,
		~g\left(X_{i}^{-}~|~X_{1:~i-1}\right)\textless s
		\text{~and~}
		Y_{i} < s
		~|~X_{1:~i-1},\xi_{1:~i-1}\right\}}\notag\\
	&=&
	\mathrm{P}\left\{X_{i}<L_{i}~|~X_{1:~i-1}\right\}
	+\mathrm{P}\left\{X_{i}=L_{i}~|~X_{1:~i-1}\right\}\notag\\
	&&\cdot\mathrm{P}\left\{\left.\xi_{i} <
	\frac{s-g\left(X_{i}^{-}~|~X_{1:~i-1}\right)}
	{g\left(X_{i}~|~X_{1:~i-1}\right)
		-g\left(X_{i}^{-}~|~X_{1:~i-1} \right)}
	~\right|~{X_{1:~i-1},X_{i}=L_{i},\xi_{1:~i-1}}
	\right\}\notag\\
	&=& g\left(L_{i}^{-}~|~X_{1:~i-1}\right)
	+\left[g\left(L_{i}~|~X_{1:~i-1}\right)
	-g\left(L_{i}^{-}~|~X_{1:~i-1}\right)\right]
	\notag \\
	&& \cdot       \frac{s-g\left(L_{i}^{-}~|~X_{1:~i-1}\right)}
	{g\left(L_{i}~|~X_{1:~i-1}\right)
		-g\left(L_{i}^{-}~|~X_{1:~i-1}\right)}
	\notag\\
	&=& s. \notag
	\end{eqnarray}
	{Note that the second equality holds since from the definition of $Y_{i}$ in Eq. (\mbox{\ref{eq. Y_{n}}}) $Y_{i}\textless s$ is equivalent to $\xi_{i}<\frac{s-g\left(X_{i}^{-}~|~X_{1:~i-1}\right)}     {g\left(X_{i}~|~X_{1:~i-1}\right)-g\left(X_{i}^{-}~|~X_{1:~i-1} \right)}$, and the latter is well defined under the condition $X_{i}=L_{i}$ suggesting that $g\left(X_{i}~|~X_{1:~i-1}\right)-g\left(X_{i}^{-}~|~X_{1:~i-1} \right)\textgreater 0$. The third equality holds since $\forall i,~\xi_{i}\sim U(0,1)$ and the series $\{\xi_i\}_{i\in\mathbb{Z}_{\textgreater 0}}$ is independent of themselves and also of the series $\{X_i\}_{i\in\mathbb{Z}_{\textgreater 0}}$, which make the conditional probability vanish, i.e., }
	\begin{eqnarray*}
		&&\mathrm{P}\left\{\left.\xi_{i} <
		\frac{s-g\left(X_{i}^{-}~|~X_{1:~i-1}\right)}
		{g\left(X_{i}~|~X_{1:~i-1}\right)
			-g\left(X_{i}^{-}~|~X_{1:~i-1} \right)}
		~\right|~X_{1:~i-1},X_{i}=L_{i},\xi_{1:~i-1}
		\right\}\\
		&=& \mathrm{P}\left\{\left.\xi_{i} <
		\frac{s-g\left(L_{i}^{-}~|~X_{1:~i-1}\right)}
		{g\left(L_{i}~|~X_{1:~i-1}\right)
			-g\left(L_{i}^{-}~|~X_{1:~i-1} \right)}
		~\right|~X_{1:~i-1},\xi_{1:~i-1}
		\right\}\\
		&=&   \frac{s-g\left(L_{i}^{-}~|~X_{1:~i-1}\right)}
		{g\left(L_{i}~|~X_{1:~i-1}\right)
			-g\left(L_{i}^{-}~|~X_{1:~i-1} \right)}.
	\end{eqnarray*}
	
	{Cases II): for a given $s\in(0,1)$ there are $g\left(L_i~|~X_{1:~i-1}\right)\geq s$ and $g\left(L_i^{-}~|~X_{1:~i-1}\right)=s$. In this case it is impossible that $ g\left(X_{i}~|~X_{1:~i-1}\right)\geq s$ while $g\left(X_{i}^{-}~|~X_{1:~i-1}\right) < s$, so for $s\in(0,1)$ Eq. (\mbox{\ref{Lemma3gao}}) changes to be}
	\begin{eqnarray*}{\label{eq estimation lemma 3 one}}
		&&\mathrm{P}\left\{Y_{i} < s~|~X_{1:~i-1},\xi_{1:~i-1}\right\} \\
		&=& \mathrm{P}\left\{
		g\left(X_{i}~|~X_{1:~i-1}\right) < s
		~|~X_{1:~i-1}\right\} \notag\\
		&=& \mathrm{P}\left\{X_{i}<L_{i}~|~X_{1:~i-1}\right\}\notag\\
		&=& g\left(L_{i}^{-1}~|~X_{1:~i-1}\right) \notag\\
		&=& s. \notag
	\end{eqnarray*}
	
	{Clearly, both cases indicate $\mathrm{P}\left\{Y_{i} < s~|~X_{1:~i-1},\xi_{1:~i-1}\right\}=s$. Note that} $$\mathrm{P}\left\{Y_{i} < s\right\}=\mathrm{E}\left[\mathrm{P}\left\{Y_{i} < s~|~X_{1:~i-1},\xi_{1:~i-1}\right\}\right]=s,$$ {we thus have $Y_{i}\sim U(0,1)$. Further, for any $l_{1},\cdots,l_{i-1},\epsilon_{1},\cdots,\epsilon_{i-1} \in \mathbb{R}$ there is}
	\begin{eqnarray}
	&&\mathrm{P}\left\{Y_{i} < s,X_{1}<l_{1},\cdots,X_{i-1}<l_{i-1},\xi_{1}<\epsilon_{1},\cdots,\xi_{i-1}<\epsilon_{i-1}\right\} \notag\\
	&=&
	\mathrm{P}\left\{Y_{i} < s|X_{1:~i-1},\xi_{1\sim i-1}\right\}\cdot\mathrm{P}\left\{X_{1}<l_{1},\cdots,X_{i-1}<l_{i-1},\xi_{1}<\epsilon_{1},\cdots,\xi_{i-1}<\epsilon_{i-1}\right\} \notag\\
	&=& \mathrm{P}\left\{Y_{i} < s\right\}
	\cdot\mathrm{P}\left\{X_{1}<l_{1},\cdots,X_{i-1}<l_{i-1},\xi_{1}<\epsilon_{1},\cdots,\xi_{i-1}<\epsilon_{i-1}\right\}. \notag
	\end{eqnarray}
	{Therefore, $\forall i\geq 1,~Y_{i}$ is independent of $X_{1},\cdots,X_{i-1},\xi_{1},\cdots,\xi_{i-1}$. A look back at Eq. (\mbox{\ref{eq. Y_{n}}}) reveals that for any $\forall k\textless i$, $Y_k$ is completely defined by $X_{1},\cdots,X_{k}$ and $\xi_{k}$, so $Y_{i}$ is independent of $Y_k$. We thus get that the series $\{Y_{i}\}_{i\in\mathbb{Z}_{\textgreater 0}}$ is mutually independent and also obey the distribution of $U(0,1)$.}
\end{proof}
\begin{proof}[Proof of Theorem 3]
	We intend to prove this theorem in two separate cases. For convenience, we will define $\liminf_{T\to+\infty} \frac{1}{T}\int_{0}^{T} \mathbbold{1}_{\left\{V(x(s)) < V_{1}\right\}} \dd s$ by $\tilde{D}$.
	
	{According to \textit{Remark \mbox{\ref{remark 4}}}, the whole evolution time can be well decomposed into the duration of each loop, as mentioned in \textit{Remark \mbox{\ref{remark 2}}}.
		Therefore, $\tilde{D}$ can be evaluated by considering the up crossing time of each loop.
	}
	
	(i) The number of loops is finite. Denote this number by $\tilde{i}+1$, then $\tau_{2\tilde{i}+2}=+\infty$.  From \textit{Lemma \ref{Lem down crossing time is finite}}, the loop down crossing time is finite, so $\tau_{2\tilde{i}+1}$ should be almost surely infinite. This implies that after $\tau_{2\tilde{i}}$ the function $V(x(t))$ is strictly less than $V_{1}$. Hence, all through the whole evolution period there is
	\begin{equation*}
	\liminf_{T\to+\infty} \frac{1}{T}\int_{0}^{T} \mathbbold{1}_{\left\{V(x(s)) < V_{1}\right\}} \dd s
	=1 \geq \frac{t_{\mathrm{u}}}{t_{\mathrm{u}}+t_{\mathrm{d}}}\quad\quad {\mathrm{a.~s.}},
	\end{equation*}
	i.e.,
	\begin{equation}{\label{eq estimation for finite many loops cases}}
	\mathrm{P}\left\{
	\left\{\text{There are \text{finitely many loops}}\right\}
	\bigcap\left\{
	\tilde{D}
	< \frac{t_{\mathrm{u}}}{t_{\mathrm{u}}+t_{\mathrm{d}}}
	\right\}
	\right\}=0.
	\end{equation}
	
	(ii) The number of loops is infinite. Denote $i(T)$ as $\max\left\{i~|~\tau_{2i}\leq T\right\}$, and $\tilde{D}$ can be estimated from definition
	\begin{eqnarray}
	\tilde{D}	  &\geq&
	\liminf_{T\to\infty}\frac{\sum_{i=0}^{i(T)-1} (\tau_{2i+1}-\tau_{2i})+\tau_{2i(T)+1}\wedge T-\tau_{2i(T)}}
	{\sum_{i=0}^{i(T)-1} (\tau_{2i+2}-\tau_{2i})+ T-\tau_{2i(T)}}   \notag \\
	&\geq&
	\liminf_{T\to\infty}
	\min\left\{
	\frac{\sum_{i=0}^{i(T)-1} (\tau_{2i+1}-\tau_{2i})}
	{\sum_{i=0}^{i(T)-1} (\tau_{2i+2}-\tau_{2i})}
	~,~
	\frac{\sum_{i=0}^{i(T)-1} (\tau_{2i+1}-\tau_{2i}) + \tau_{2i(T)+1} -\tau_{2i(T)}}
	{\sum_{i=0}^{i(T)-1} (\tau_{2i+2}-\tau_{2i})+ \tau_{2i(T)+2}-\tau_{2i(T)}}
	\right\}  \notag \\
	&=&  \liminf_{n\to\infty}
	\frac{\sum_{i=0}^{n}(\tau_{2i+1}-\tau_{2i})}
	{\sum_{i=0}^{n}(\tau_{2i+1}-\tau_{2i})
		+\sum_{i=0}^{n}(\tau_{2i+2}-\tau_{2i+1})} \notag \\
	&=& \liminf_{n\to\infty}
	\frac{1}{1+
		\frac{\sum_{i=0}^{n}(\tau_{2i+2}-\tau_{2i+1})/(n+1)\cdot(n+1)/n}
		{\sum_{i=1}^{n}(\tau_{2i+1}-\tau_{2i})/n+(\tau_{1}-\tau_{0})/n}
	} \notag \\
	&\geq&
	\frac{1}{1+
		\frac{\limsup_{n\to\infty}\sum_{i=0}^{n}(\tau_{2i+2}-\tau_{2i+1})/(n+1)}
		{\liminf_{n\to\infty}\sum_{i=1}^{n}(\tau_{2i+1}-\tau_{2i})/n}
	}, \notag
	\end{eqnarray}
which means if
$\tilde{D}< \frac{t_{\mathrm{u}}}{t_{\mathrm{u}}+t_{\mathrm{d}}}$ happens, then $\frac{1}{1+
		\frac{\limsup_{n\to\infty}\sum_{i=0}^{n}(\tau_{2i+2}-\tau_{2i+1})/(n+1)}
		{\liminf_{n\to\infty}\sum_{i=1}^{n}(\tau_{2i+1}-\tau_{2i})/n}}< \frac{t_{\mathrm{u}}}{t_{\mathrm{u}}+t_{\mathrm{d}}}$ must happen, i.e., $$\frac{\limsup_{n\to\infty}\sum_{i=0}^{n}(\tau_{2i+2}-\tau_{2i+1})/(n+1)}
		{\liminf_{n\to\infty}\sum_{i=1}^{n}(\tau_{2i+1}-\tau_{2i})/n}>\frac{t_{\mathrm{d}}}{t_\mathrm{u}}$$  must happen. We thus have the relationship between their occurrence events
	\begin{eqnarray}
	\left\{\tilde{D}< \frac{t_{\mathrm{u}}}{t_{\mathrm{u}}+t_{\mathrm{d}}}\right\}\notag 
	&\subseteq&\left\{\frac{\limsup_{n\to\infty}\sum_{i=0}^{n}(\tau_{2i+2}-\tau_{2i+1})/(n+1)}
	{\liminf_{n\to\infty}\sum_{i=1}^{n}(\tau_{2i+1}-\tau_{2i})/n}\textgreater \frac{t_{\mathrm{d}}}{t_\mathrm{u}}	\right\}\notag \\
	&\subseteq&\left\{\liminf_{n\to\infty} \sum_{i=1}^{n} \frac{\tau_{2i+1}-\tau_{2i}}{n} {\textless} t_{\mathrm{u}}\right\}\bigcup\left\{\limsup_{n\to\infty} \sum_{i=0}^{n} \frac{\tau_{2i+2}-\tau_{2i+1}}{n+1}{\textgreater} t_{\mathrm{d}}\right\}\notag\\
	&\subseteq&\left\{\liminf_{n\to\infty} \sum_{i=1}^{n} \frac{\tau_{2i+1}-\tau_{2i}}{n} {\textless} t_{\mathrm{u}}\right\}\bigcup\left\{\limsup_{n\to\infty} \sum_{i=1}^{n} \frac{\tau_{2i}-\tau_{2i-1}}{n} {\textgreater} t_{\mathrm{d}}\right\}.\notag
	\end{eqnarray}
	{Then we have}
	\begin{eqnarray}\label{eq estimation for infinite many loops cases}
	&&\mathrm{P}\left\{
	\left\{\text{There are \text{infinitely many loops}}\right\}
	\bigcap\left\{
	\tilde{D}
	< \frac{t_{\mathrm{u}}}{t_{\mathrm{u}}+t_{\mathrm{d}}}
	\right\}
	\right\} \\
	&\leq& 	\mathrm{P}\left\{
	\left\{\text{\textup{There are \text{infinitely many loops}}}\right\}
	\bigcap\left\{
	\liminf_{n\to\infty} \sum_{i=1}^{n} \frac{\tau_{2i+1}-\tau_{2i}}{n} {\textless} t_{\mathrm{u}}
	\right\}
	\right\} \notag\\
	&&+
	\mathrm{P}\left\{
	\left\{\text{\textup{There are \text{infinitely many loops}}}\right\}
	\bigcap\left\{
	\limsup_{n\to\infty} \sum_{i=1}^{n} \frac{\tau_{2i}-\tau_{2i-1}}{n} {\textgreater} t_{\mathrm{d}}
	\right\}
	\right\}\notag \\
	&=&0, \notag
	\end{eqnarray}
	{where the last equality holds by applying the results of \textit{Propositions} \mbox{\ref{Lem up crossing time is almost surely lower bounded}} and \mbox{\ref{Lem down crossing time is almost surely upper bounded}}. Namely, in case ii) $\tilde{D}\geq\frac{t_{\mathrm{u}}}{t_{\mathrm{u}}+t_{\mathrm{d}}}$ is almost surely.}
	
	Finally, combining Eqs. (\ref{eq estimation for finite many loops cases}) and (\ref{eq estimation for infinite many loops cases}) we prove the result true.
\end{proof}

\end{document}